\documentclass[absolute]{mymathart}
\usepackage[theorems]{mymathmacros}
\usepackage[usenames,dvipsnames]{color}
\usepackage{subfig}
\usepackage[shortlabels]{enumitem}
\usepackage{hyperref}
\usepackage{amsmath}
\usepackage{amssymb}
\usepackage{graphicx}
\usepackage{epigraph}
\usepackage{amsfonts,amssymb,color}
\usepackage[mathscr]{eucal}
\usepackage{amsmath, amsthm}
\usepackage{mathrsfs}
\usepackage{amsbsy}
\usepackage{float}
\usepackage{verbatim}
\usepackage{amsfonts} 
\usepackage{pifont}
\usepackage{tikz-cd}
\usepackage{color}

\title[On the finite factorial numbers]{On the finite factorial numbers}

\author{Liangang Ma}

\address{Liangang Ma, School of Mathematics and Statistics, Ludong University, Yantai 264025, Shandong, P. R. China.} 
\email{maliangang000@163.com}

\thanks{The work is supported by NSFC-12001056 and LU-20220028.}  

\newtheorem{theorem}[subsection]{Theorem}
\newtheorem{lemma}[subsection]{Lemma}
\newtheorem{pro}[subsection]{Proposition}

\newtheorem*{Bird's Theorem}{Bird's Theorem}
\newtheorem*{Evertse's Theorem}{Evertse's Theorem}
\newtheorem*{Schmidt Subspace Theorem}{Schmidt Subspace Theorem}

\newtheorem{coro}[subsection]{Corollary}
\newtheorem{exa}[subsection]{Example}
\newtheorem{rem}[subsection]{Remark}

\counterwithout{footnote}{section}

\numberwithin{equation}{section}

\begin{document} 

\begin{abstract}
We introduce two families of transcendental numbers which we call finite factorial (FF) and partially finite factorial (PFF) numbers respectively, with the former one being subfamily of the latter one. These numbers arise naturally from some transcendental criterion for real numbers via their $b$-ary expansions. We show that rational numbers (eventually periodic words) can not be finite factorial. Then we consider the geometric (topological) properties of the collection of all the FF numbers, including its countability, density and Hausdorff dimension. Some numerical examples are given to illustrate certain results in the work.           
\end{abstract}
 
 \maketitle

\section{Introduction}\label{sec1}

We focus on reading \emph{transcendence\textbackslash algebraicity} of real numbers from their $b$-ary expansions in this work. The topic dates back to Liouville (\cite{Lio}), who proved that the decimal number
$$(.1100010000000000000000001000\cdots)_{10}=\sum_{i=1}^\infty 10^{-i!}$$
is transcendental. The famous Liouville's inequality which takes core role in his proof is then consistently polished by  Thue (\cite{Thu}), Siegel (\cite{Sie}), Dyson (\cite{Dys}) and Roth (\cite{Rot}). The Schmidt Subspace Theorem (\cite{Schm}) can be understood along this line, with its application in bounding numbers of solutions of certain Diophantine equations. 

Although almost all real numbers are transcendental, it is not an easy work to decide the transcendence or algebraicity of an individual one.  It is due to Borel (\cite{Bor}) that almost all real numbers are \emph{normal} (numbers such that any length-$k$ word (sequence) in their $b$-ary expansions appears with probability $1/b^k$). It is a big conjecture that all irrational algebraic numbers are normal, which offers inspiration for many transcendence criteria. Ferenczi and Mauduit (\cite{FM}) proved that irrational numbers whose $b$-ary expansions admit some initial \emph{powers} (refer to \cite{BHZ} by Berth\'e-Holton-Zamboni) are transcendental. There are generalizations of their criterion, for example, by     
Allouche-Zamboni (\cite{AZ}) and Adamczewski-Bugeaud-Luca (\cite{ABL}), in various contexts. 

It is Cobham (\cite{Cob}) who first related the $b$-ary expansions of irrational numbers to the particular Turing machines: the finite \emph{automata}. One can refer to Allouche and Shallit's book \cite{AS1} for a systematic account of the topic. The Cobham-Loxton-van der Poorten conjecture (\cite{LV1, LV2}) asserts that an irrational algebraic
number cannot be generated by a finite automaton, which is finally obtained by Adamczewski and Bugeaud (refer to \cite{AB1}, see also \cite{AF} by Adamczewski and Faverjon for an alternative proof). Note that some special automatic numbers are known to be transcendental for a long time, for example, numbers with binary expansions being the Thue-Morse word (refer to Mahler \cite{Mah}), the Fibonacci word and the  Rudin–Shapiro word (refer to Allouche-Shallit \cite{AS2}). 

The distribution of finite words in an infinite automatic word is an independent topic from some other points of views.     
Cobham (\cite{Cob}) showed that the density of any letter in an automatic word is necessarily rational, in case of its existence. He also showed that the logarithmic density of letters in automatic words always exists. This is extended to all \emph{morphic} words by Bell (\cite{Bel1}), which answers a question of Allouche and Shallit. Bell also showed that any $k$-automatic set admits rational  upper and lower density regarding as a subset of the natural numbers (\cite{Bel2}). There are some profound programme posed by Gelfond (\cite{Gel}) on the density of letters along polynomial subsequences of certain sequences ($mod$ some integer). Complete answers are given by Mauduit and Rivat along linear and quadratic subsequences (\cite{MR1, MR2}), while there are only some partial results along higher-degree subsequences by  Drmota, Mauduit, and Rivat, assuming sufficient largeness of the base and coprimality of the base with leading coefficient of the polynomial (\cite{DMR}). Recently Adamczewski, Drmota and M\"ullner discover that the existence of (logarithmic) density along regularly varying subsequences depends on the property of the \emph{final strongly connected components} of the automata (\cite{ADM}).

In the following we introduce a new transcendence criterion from the $b$-ary expansions of real numbers. Classical criteria typically focus on the denominators of the rational truncations from the $b$-ary expansions of numbers, while in this work we lay emphasis on both the numerators and denominators of the rational truncations simultaneously (refer to  \cite{AB3} for another example on exploiting the subtle interaction between them). Without loss of generality we restrict our attention to numbers in the unit interval $(0,1)$.

\begin{theorem}\label{thm1}
For an integer base $b\geq 2$ and an irrational real number $0<\xi<1$ with its $b$-ary expansion $\bold{a}=a_1a_2a_3\cdots$, if there exist a finite set $S=\{s_1, s_2, \cdots, s_l\}$ of primes and a strictly increasing sequence $\{n_i\}_{i=1}^\infty$ of positive integers such that 
\begin{equation}\label{eq5}
\cfrac{c_{u_{n_i},S}}{u_{n_i}}\geq b^{(\epsilon_{1.1}-1)n_i}
\end{equation}
for any $i$ and some positive number $0<\epsilon_{1.1}\leq 1$,  in which $u_{n_i}$ is the numerator of the $n_i$-th rational truncation of $\xi$ and $c_{u_{n_i},S}$ is the $S$-component of $u_{n_i}$,  then $\xi$ is a transcendental number.
\end{theorem}

The proof of the transcendence criterion relies on the $p$-adic Schmidt Subspace Theorem due to Schlickewei (\cite{Schl}), refer to the Appendix. Theorem \ref{thm1} demonstrates importance of prime factorization of (sequences of) integers, while one can find some interesting results on prime factorization of (sequences of) sums of integers in \cite{GST} by Gy\"ory, Stewart and Tijdeman. 

In the following we call a real (rational and irrational) number satisfying (\ref{eq5})  a \emph{partially finite factorial} (PFF) number with respect to some $b, S, \{n_i\}_{i=1}^\infty, \epsilon_{1.1}$. In case (\ref{eq5}) holds for $\epsilon_{1.1}=1$ with respect to some $b, S, \{n_i\}_{i=1}^\infty$ it is named a \emph{finite factorial} (FF) number. Theorem 
\ref{thm1} demonstrates its power in deciding the transcendence of reals in case one collects plenty information about the prime factorization of the numerators of the rational truncations of the reals. 

\begin{theorem}\label{thm5}
Let $b\geq 2$ be some integer base with its prime factorization $b=s_1^{k_1}s_2^{k_2}\cdots s_l^{k_l}$ for some finite set of primes $S=\{s_1, s_2, \cdots, s_l\}$ with positive integers $\{k_1,k_2,\cdots,k_l,l\}$. Let 
\begin{center}
$\{m_i=s_1^{r_{i,1}}s_2^{r_{i,2}}\cdots s_l^{r_{i,l}}\}_{i=1}^\infty$
\end{center}
be a sequence of nodes with 
\begin{equation}\label{eq27}
\inf_{1\leq i<\infty}\Big\{\cfrac{\sum_{j=1}^l r_{i+1,j}}{\sum_{j=1}^l r_{i,j}}\Big\}\geq 1+\epsilon_{1.2}
\end{equation}
for some fixed $\epsilon_{1.2}>0$.  Then the number with $b$-ary expansion
\begin{center}
$\bold{a}_1\bold{a}_2\bold{a}_3\cdots$
\end{center}
in $(0,1)$ is transcendental, in which $\bold{a}_i$ is the $b$-ary expansion of the node $m_i$ for any $1\leq i<\infty$.
\end{theorem}

Theorem \ref{thm5} covers some classical results, while some numbers with completely new expansions can be added to the collection of transcendentals in virtue of it. For example, let $b=10=2\cdot5$ and $m_i=2^{(i+1)!-i!-1}5^{(i+1)!-i!-1}$ for any $1\leq i<\infty$. Note that the $10$-ary expansion of $m_1$ is $1$, the $10$-ary expansion of $m_2$ is $1000$, the $10$-ary expansion of $m_3$ is $1000000000000000000$, $\cdots$. Then the Liouville number
$$(.1100010000000000000000001000\cdots)_{10}=\sum_{i=1}^\infty 10^{-i!}$$
can be concluded to be transcendental via Theorem \ref{thm5}. However, the core contribution of Theorem \ref{thm5} lies in its inconsistent cases. For example, let $b=10=2\cdot5$ and $m_i=2^{(i+1)!-i!-1}$ for any $1\leq i<\infty$, then the number 
$$(.1813107239614081257132168796771975168...)_{10}$$
is transcendental, in virtue of Theorem \ref{thm5}. More interesting examples are provided in Section \ref{sec9}.

Although it is easy to distinguish rational numbers from the irrational ones via their $b$-ary expansions, it is an interesting question to ask whether rational numbers are finite factorial or partially finite factorial. In virtue of a theorem of Evertse (see Evertse's Theorem in the Appendix) we prove the following result.

\begin{theorem}\label{pro3}
A rational number (or an eventually periodic word in $\bold{A}_{b}$) can never be finite factorial, with respect to any base $b\geq 2$. 
\end{theorem}

Elementary proof does apply to some special case (refer to Proposition \ref{pro4}) of  Theorem \ref{pro3}, while we can not prove Theorem \ref{pro3} without Evertse's Theorem. We conjecture that rational numbers can not be partially finite factorial, however, this will not be dealt with in this work. 

Let $\Xi^{FF}$ and $\Xi^{PFF}$ be the collections of all the finite factorial and partially finite factorial numbers in the unit interval respectively, so $\Xi^{FF}\subset \Xi^{PFF}$.  In the following we turn to the geometric (topological) properties of the two sets $\Xi^{FF}$ and $\Xi^{PFF}$. As to the finite factorial numbers, we have the following result.

\begin{theorem}\label{thm3}
The set $\Xi^{FF}$ is an uncountable set of null Hausdorff dimension, which is dense in $[0,1]$.
\end{theorem}

The Bird's Theorem (refer to Bird's Theorem in the Appendix) is a critical instrument in establishing the uncountability as well as the density of the set. This means we can approximate any real number by sequences of finite factorial numbers. The result may induce some interesting way in approximation of (algebraic or transcendental) numbers. Traditional point of view focuses on approximating transcendental numbers by certain algebraic ones, for example, one can refer to Davenport-Schmidt (\cite{DS}) and Roy-Waldschmidt (\cite{RW1, RW2}). 

Good's technique (\cite{Goo}) on decomposing certain sets arising from continued fraction expansions fits well into our estimation on the upper bound of the Hausdorff dimension of the set.

As a superset of $\Xi^{FF}$, the set of partially finite factorial numbers $\Xi^{PFF}$ inherits its uncountability and density spontaneously.

\begin{coro}\label{thm4}
The set $\Xi^{PFF}$ is an uncountably dense subset of $[0,1]$.
\end{coro}  

Further properties of the set of partially finite factorial numbers will be developed in a successive work.  

The layout of the work is as following. We introduce some concepts and notations on operations of integers, $b$-ary expansions of numbers and Hausdorff dimension in Section  \ref{sec8}, while some other local concepts and notations are given intermittently in following sections. In Section \ref{sec2} we justify the transcendence criterion: Theorem \ref{thm1}. In Section \ref{sec9} we prove Theorem \ref{thm5}, with some ameliorative versions of it presented.  The readers will see more interesting explicit examples of transcendentals built in the section. In Section \ref{sec7} we clarify the relationship between rational numbers and the FF numbers by Theorem \ref{pro3}. The following section is dedicated to the geometric (topological) properties of the FF numbers, Theorem \ref{thm3} is proved here. In Section \ref{sec5} we provide the readers with some examples on the \emph{ancestor-descendant} relationship (\emph{cf.} Section \ref{sec3}) and the least FF number beginning with some given word (\emph{cf.} Proposition \ref{thm2}). The Appendix section contains statements of (some version of) the Schmidt Subspace Theorem, the Evertse's Theorem and the Bird's Theorem, which will be employed to tackle relevant problems in due course.

\section{Concepts and notations}\label{sec8}

Let $\mathbb{N}=\{0,1,2,\cdots\}, \mathbb{Z}, \mathbb{Q}, \mathbb{R}$ be the set of natural numbers, integers, rationals and real numbers respectively, let $*_+=*\cap (0,\infty)$ be the set of positive members for any $*=\mathbb{N}, \mathbb{Z}, \mathbb{Q}, \mathbb{R}$. We denote by $p_i$ the $i$-th prime number with $p_1=2$ for any $i\in\mathbb{N}_+$. The collection of initial $k$ primes is denoted by $P_k=\{p_i\}_{i=1}^k$, with the collection of all primes being $P=\cup_{k=1}^\infty P_k$.   For an integer $r\geq 2$, let $F_r$ be its collection of prime factors, that is, 
$$F_r=\{s: s|r \text{ is a prime}\},$$ 
in which $|$ means dividing. The greatest common divisor of two positive integers $r_1$ and $r_2$ is denoted by $gcd(r_1,r_2)\in\mathbb{N}_+$. We also write $r_1\perp r_2$ if $gcd(r_1,r_2)=1$, that is, $r_1$ is coprime with $r_2$.

For some prime $s$ and some rational $\cfrac{r_1}{r_2}$ with $r_1, r_2\in\mathbb{N}_+$, let 
\begin{center}
$|\cfrac{r_1}{r_2}|_s=\min\{s^{-j}: s^{j}|r_1, j\in\mathbb{N}\}\cdot\max\{s^j: s^{j}|r_2, j\in\mathbb{N}\}$
\end{center} 
be the $s$-\emph{adic absolute value} of $\cfrac{r_1}{r_2}$ (note that it is possible $|\cfrac{r_1}{r_2}|_s=1$). 

We traditionally use 
$$S=\{s_1, s_2, s_3, \cdots, s_l\}\subset P$$ 
to denote a finite set of $l=\#S$ primes in strictly increasing order with respect to their subscripts, which means $2 \leq s_1<s_2<s_3<\cdots<s_l$.  

 For a finite set $S=\{s_1, s_2, \cdots, s_l\}$ of primes and  a positive integer $r$, the $S$-\emph{component} of $r$ is defined to be the integer 
$$c_{r,S}=\max\{\prod_{j=1}^l s_j^{i_j}: \prod_{j=1}^l s_j^{i_j}|r,\ (i_1, i_2, \cdots, i_l)\in\mathbb{N}^l\}.$$

We traditionally employ lower(upper)-case letters to denote some integer (set of integers) and lower-case Greece letters to denote arbitrary real numbers. The bold letters are customarily employed to denote (finite or infinite) words (or sequences) which are concatenations of letters.  For a real number $\alpha\in\mathbb{R}$, let $\lfloor\alpha\rfloor$ be its integer part and $(\alpha)=\alpha-\lfloor\alpha\rfloor$ be its fractional part. The fractional part decides the  transcendence\textbackslash algebraicity of a real number.

In the following we introduce some notations on $b$-ary expansions of real numbers. Since we focus on the algebraicity and transcendence of real numbers in the unit interval $(0,1)$ in this work, we restrict to the $b$-ary expansions of numbers in $(0,1)$. Throughout the work the base $b$ will mean some integer greater than or equal to $2$. Let $B=\{0,1,2,\cdots,b-1\}$ be the \emph{alphabet} for some base $b\geq 2$, with its members named \emph{letters}. Let $\bold{A}_b^F=\cup_{i\in\mathbb{N}} B^i$ be the collection of all finite words with respect to the base $b$, in which $B^0=\emptyset$ containing only the empty word. Let 
\begin{center}
$
\begin{array}{ll}
& \bold{A}_{b,0}\\
=&\{\bold{a}=a_1a_2a_3\cdots\in B^\infty: \mbox{ there exists some\ } i\in\mathbb{N}_+, \mbox{such that\ } a_j=0 \mbox{\ for any\ } j\geq i\}
\end{array}
$
\end{center}
be the collection of infinite words which are concatenations of finite words in $\bold{A}_{b}^F$ with the infinite word $\bold{0}=000\cdots$ ($\bold{k}=kkk\cdots$ for some letter $k\in B$ may be a finite or infinite word, while one can judge its length easily in due contexts) for any base $b\geq 2$. Let 
$$\bold{A}_{b}=B^\infty\setminus \bold{A}_{b,0}.$$ 
The \emph{length} $|\bold{a}|$ of a (finite or infinite) word $\bold{a}\in \bold{A}_{b}\cup\bold{A}_{b}^F$ is the number of letters in the word (it is $\infty$ for an infinite word).   For an infinite (finite) word $\bold{a}=a_1a_2a_3\cdots\in \bold{A}_{b}$ ($\bold{a}\in \bold{A}_{b}^F$), let 
$$\bold{a}|_k=a_1a_2a_3\cdots a_k\in\bold{A}_{b}^F$$
be its $k$-th restriction  for $1\leq k< \infty$ ($1\leq k\leq |\bold{a}|$).

For some base $b$ and some real number $\xi\in (0,1)$, the $b$-ary expansion of $\xi$ is the unique infinite word 
$$\bold{a}=\bold{a}_\xi=a_1a_2a_3\cdots$$ 
in $\bold{A}_{b}$ such that 
$$\xi=\cfrac{a_1}{b}+\cfrac{a_2}{b^2}+\cfrac{a_3}{b^3}+\cdots=(.a_1a_2a_3\cdots)_b.$$

By ruling out words in $\bold{A}_{b,0}$ we guarantee the uniqueness of the expansions. Alternatively, some people prefer to rule out concatenations of finite words with the infinite word $\bold{b-1}=b-1b-1b-1\cdots$, which imposes merely some little distinct interpretation on some results below.

The concept applies to integers with some mild distinction. For an integer $u\in\mathbb{N}_+$, the $b$-ary expansion of $u$ is the unique finite word 
$$\bold{a}=\bold{a}_u=a_1a_2a_3\cdots a_k$$ 
in $B^k$ with $a_1\neq 0$ such that 
$$u=u_{\bold{a}}=a_1b^{k-1}+a_2b^{k-2}+a_3b^{k-3}+\cdots+a_k=(a_1a_2a_3\cdots a_k.)_b$$
for some $k\in\mathbb{N}_+$. 

We will frequently omit the symbols $()_b$ in case the base being obvious in the following. Conversely, given any infinite (finite) word $\bold{a}\in \bold{A}_{b}$ ($\bold{a}\in B^k$) for some base $b$ (and $k\in\mathbb{N}_+$), there is a unique number $\xi_{\bold{a}}\in (0,1)$ ($u_{\bold{a}}\in\mathbb{N}_+$) whose $b$-ary expansion is $\bold{a}$. There is some relationship between the $b$-ary expansions of reals in the unit interval and $\mathbb{N}_+$. For a real number $\xi_{\bold{a}}=(.\bold{a})_b=(.a_1a_2a_3\cdots)_b$ with respect to some base $b$, considering the restriction $\bold{a}|_k$ for some $k\in\mathbb{N}_+$, the integer
$$u_{\bold{a}|_k}=u_{k,\bold{a}}=u_{k,\xi}=\sum_{i=1}^k a_i b^{k-i}$$
is the numerator of its $k$-th rational truncation 
$$(.a_1a_2a_3\cdots a_k)_b=(.a_1a_2a_3\cdots a_k000\cdots)_b=\cfrac{u_{\bold{a}|_k}}{b^k},$$
so we have
$$\xi_{\bold{a}}=\lim_{k\rightarrow\infty} \cfrac{u_{\bold{a}|_k}}{b^k}$$
in fact. We will simply write $u_{\bold{a}|_k}=u_k$, omitting the symbol $\bold{a}$ or $\xi$, in case the infinite word $\bold{a}$ or number $\xi$ being clear in the following. 

It is well-known that rational numbers admit eventually periodic $b$-ary expansions for any base $b\geq 2$ (see for example \cite{HW} by Hardy-Wright).

In order to deal with the Hausdorff dimension of the set of FF numbers, we introduce some related concepts and notations.  For a set $\Xi\subset (0,1)$, let $|\Xi|$ be its diameter with respect to the Euclidean metric. 
For $\alpha,\delta\in\mathbb{R}_+$, let 
$$H^\alpha_\delta(\Xi)=\inf\big\{\sum_{i=1}^\infty |\Xi_i|^\alpha: \{\Xi_i\}_{i\in\mathbb{N}_+} \mbox{ is a } \delta\mbox{-cover of\ } \Xi\big\}$$ 
and 
$$H^\alpha(\Xi)=\lim_{\delta\rightarrow 0} H^\alpha_\delta(\Xi)$$
be the $\alpha$ dimensional \emph{Hausdorff measure} of $\Xi$. The \emph{Hausdorff dimension} $dim_H \Xi$ (as an indicator for the size and complexity of $\Xi$, especially when $\Xi$ is of null measure) is defined as
$$dim_H \Xi=\inf\{\alpha: H^\alpha(\Xi)=\infty\}=\sup\{\alpha: H^\alpha(\Xi)<\infty\}.$$
For any finite word $\bold{a}_k=a_1a_2\cdots a_k\in B^k$ of length $k\in\mathbb{N}_+$, let 
$$I_{\bold{a}_k}=\{\xi=(.a_{1,\xi}a_{2,\xi}\cdots a_{k,\xi} a_{k+1,\xi}\cdots)_b\in(0,1): a_{1,\xi}a_{2,\xi}\cdots a_{k,\xi}=\bold{a}_k\}$$
be the corresponding depth-$k$ closed fundamental interval with respect to some base $b$ (except the two ones with terminals $0,1$). Let 
$$I_{\bold{a}_k}^o=I_{\bold{a}_k}\setminus\{.\bold{a}_k 000\cdots, .\bold{a}_k (b-1)(b-1)(b-1)\cdots\}$$
be the corresponding open interval. The fundamental intervals provide natural covers for sets arising from $b$-ary expansions.

\section{A transcendental criterion via the $b$-ary numerical system}\label{sec2}

We start from proving Theorem \ref{thm1} in this section.

\begin{proof}[Proof of Theorem \ref{thm1}]
We achieve our goal by reduction to absurdity, so we assume $\xi$ is an algebraic number from now on. Without loss of generality we can suppose the sequence $\{n_i\}_{i=1}^\infty$ satisfies 
\begin{equation}\label{eq3}
\max\{a_i: 1\leq i\leq \lfloor \epsilon_{1.1}' n_i\rfloor\}>0
\end{equation}
for some $\epsilon_{1.1}'<\epsilon_{1.1}$ and any $i\in\mathbb{N}_+$, or else we can find a subsequence with the above property and the proof follows.   By truncating the $b$-ary expansion of $\xi$ at $n_i$, we have
\begin{equation}\label{eq1}
0<\xi-\cfrac{u_{n_i}}{b^{n_i}}<\cfrac{1}{b^{n_i}}
\end{equation}
for any $i\in\mathbb{N}_+$. Considering the adic absolute values of  $u_{n_i}$ and $b^{n_i}$, we have
\begin{equation}\label{eq2}
\prod_{s\in S\cup F_b} |u_{n_i}b^{n_i}|_s=\cfrac{1}{c_{u_{n_i},S\cup F_b}b^{n_i}},
\end{equation}
in which $c_{u_{n_i},S\cup F_b}$ is the $(S\cup F_b)$-component of the integer $u_{n_i}$. Now (\ref{eq1}) and (\ref{eq2}) together give
\begin{center}
$(\prod_{s\in S\cup F_b} |u_{n_i}b^{n_i}|_s) b^{n_i}|b^{n_i}\xi-u_{n_i}|<\cfrac{1}{c_{u_{n_i},S\cup F_b}}$
\end{center}
for any $i\in\mathbb{N}_+$. By the assumption (\ref{eq3}) we have
\begin{center}
$\cfrac{1}{u_{n_i}}\leq \cfrac{1}{b^{(1-\epsilon_{1.1}')n_i}}$.
\end{center}
Considering (\ref{eq5}), it follows then 
\begin{equation}\label{eq4}
(\prod_{s\in S\cup F_b} |u_{n_i} b^{n_i}|_s) b^{n_i}|b^{n_i}\xi-u_{n_i}|<\cfrac{1}{b^{(\epsilon_{1.1}-\epsilon_{1.1}')n_i}}
\end{equation}
for any $i\in\mathbb{N}_+$. Note that $\epsilon_{1.1}-\epsilon_{1.1}'>0$. Now we consider the following two linear forms
\begin{center}
$
\begin{array}{ll}
\mathcal{L}_1(x_1,x_2)=x_1, \\
\mathcal{L}_2(x_1,x_2)=\xi x_1-x_2 
\end{array}
$
\end{center}
which are linearly independent of algebraic coefficients. In case $x_1=b^{n_i}$ and $x_2=u_{n_i}$, (\ref{eq4}) gives
\begin{center}
$(\prod_{s\in S\cup F_b} |x_1|_s|x_2|_s) |\mathcal{L}_1(x_1,x_2)||\mathcal{L}_2(x_1,x_2)|<\cfrac{1}{\max\{x_1,x_2\}^{\epsilon_{1.1}-\epsilon_{1.1}'}}$
\end{center}
for any $i\in\mathbb{N}_+$. Applying the Schmidt Subspace Theorem (see Schmidt Subspace Theorem in the Appendix) in our case, the tuple $(b^{n_i}, u_{n_i})$ must lie in finitely many proper subspaces of $\mathbb{Q}^2$. This means there exist two non-zero integers $r_1, r_2$ and an infinite subsequence $\{n_{i_j}\}_{j=1}^\infty\subset\{n_i\}_{i=1}^\infty$ such that  
\begin{equation}
r_1b^{n_{i_j}}+r_2u_{n_{i_j}}=0
\end{equation}
for any $j\in\mathbb{N}_+$. By letting $j\rightarrow\infty$, we achieve that 
$$\xi=\lim_{j\rightarrow\infty}\cfrac{u_{n_{i_j}}}{b^{n_{i_j}}}=-\cfrac{r_1}{r_2}$$
is a rational number, which contradicts the irrationality of $\xi$. 
\end{proof}

Theorem \ref{thm1} induces the following result instantly.

\begin{coro}\label{cor1}
For a real irrational number $0<\xi<1$ with its b-ary expansion $\bold{a}=a_1a_2a_3\cdots$, if there exists an infinite strictly increasing sequence $\{n_i\}_{i=1}^\infty$ of positive integers such that
$$\#\cup_{i=1}^\infty F_{u_{n_i}}<\infty,$$
then $\xi$ is transcendental.
\end{coro}

\begin{defn}
An infinite word  in $\bold{A}_{b}$ admitting a sequence $\{n_i\}_{i=1}^\infty$ with 
$$\#\cup_{i=1}^\infty F_{u_{n_i}}<\infty$$
is said to satisfy the finite factorial property (with respect to some base $b\geq 2$). It is said to satisfy the partially finite factorial property if there exist some finite set $S=\{s_1, s_2, \cdots, s_l\}$, some sequence $\{n_i\}_{i=1}^\infty$ of strictly increasing integers and some $0<\epsilon_{1.1}\leq 1$ such that (\ref{eq5}) holds (with respect to some base $b\geq 2$).
\end{defn}

These definitions apply to numbers.  A number is called a (partially) finite factorial number if it admits some (partially) finite factorial expansion with respect to some base $b\geq 2$.  It would be interesting to try to decide whether some classical words, for example, the Fibonacci word and the Thue-Morse word, are (partially) finite factorial or not.  For any  base $b\geq 2$, the Liouville word is 
$$11000100000000000000000100\cdots,$$
in which any letter is $0$ except the $i!$-th one being $1$ for any $k\in\mathbb{N}_+$.

\begin{pro}
The Liouville word is partially finite factorial with respect to any base $b\geq 2$. 
\end{pro}
\begin{proof}
Let $S=F_b$ for any base $b\geq 2$. Consider the numerators $\{u_{n_i}\}_{i=1}^\infty$ with respect to the subscripts $\{n_i=(i+3)!-1\}_{i=1}^\infty$. It is easy to see that $b\perp u_{n_i+1}$ for any $i\in\mathbb{N}_+$, so
$$c_{u_{n_i},S}=b^{(i+3)!-(i+2)!-1}$$ 
for any $i\in\mathbb{N}_+$. Then 

$$\cfrac{c_{u_{n_i},S}}{u_{n_i}}>\cfrac{b^{(i+3)!-(i+2)!-1}}{b^{(i+3)!-1}}=b^{-(i+2)!}> b^{-\frac{2}{i+3}((i+3)!-1)}=b^{(1-\frac{2}{i+3}-1)((i+3)!-1)}\geq b^{(\frac{1}{2}-1)n_i}$$
for any $i\geq 1$. Note that $\{1-\cfrac{2}{i+3}\}_{i=1}^\infty$ is a strictly increasing sequence tending to $1$ as $i\rightarrow\infty$, so (\ref{eq5}) is satisfied for the sequence $\{(i+3)!-1\}_{i=1}^\infty$ with $\epsilon_{1.1}=\cfrac{1}{2}$.
\end{proof}

We can guarantee (\ref{eq5})holds for some $\epsilon_{1.2}$ arbitrarily close to (but not equal) $1$ by considering sequences  
$\{n_i=(i+j)!-1\}_{i=1}^\infty$ for $j$ large enough.
\begin{rem}
We conjecture that the Liouville words can not be finite factorial with respect to any base $b\geq 0$. 
\end{rem}

It would be interesting to consider transplanting transcendence\textbackslash algebraicity criteria  between different numeration systems. For example, Adamczewski and Bugeaud (\cite{AB2}) showed that if the continued fraction expansion of a real number begins with arbitrarily large \emph{palindromes}, then it must be transcendental. They conjecture that one can still obtain transcendence of numbers from their palindromic pattern of $b$-ary expansions, see also Fischler \cite{Fis1, Fis2}. We wonder whether some counterpart transcendence criteria as Theorem \ref{thm1} exist in, for example, $\beta$-expansions or continued fraction expansions. 

Since (P)FF transcendental numbers obey some simple structure associated with their $b$-ary expansions, another interesting question is the relationship between (P)FF transcendental numbers and finite automaton (morphism): can (P)FF transcendental numbers be generated by finite automata (or morphisms)?

\section{Applications of the criterion to reals with explicit expansions}\label{sec9}

In this section we apply our transcendental criterion to some specific numbers expressed in the $b$-ary expansions to obtain their transcendence. Various examples are presented via Theorem \ref{thm5} as well as its ameliorative versions.  We aim to justify Theorem \ref{thm5} in this section. To do this we first establish some preliminary results. 

We borrow the notations from Theorem \ref{thm5} from now on, that is,  $b=s_1^{k_1}s_2^{k_2}\cdots s_l^{k_l}$ is some base for some finite set of primes $S=\{s_1, s_2, \cdots, s_l\}$ with positive integers $\{k_1,k_2,\cdots,k_l,l\}$ while $\{m_i=s_1^{r_{i,1}}s_2^{r_{i,2}}\cdots s_l^{r_{i,l}}\}_{i=1}^\infty$
is a sequence of integers satisfying 
$$
\liminf_{1\leq i<\infty}\{\cfrac{\sum_{j=1}^l r_{i+1,j}}{\sum_{j=1}^l r_{i,j}}\}\geq 1+\epsilon_{1.2}
$$
for some fixed $\epsilon_{1.2}>0$. We recall here that the elements of $S$ are of strictly increasing order with respect to their subscripts. To simplify some symbols let 
$$d_i=\sum_{j=1}^l r_{i,j}$$
for any $1\leq i<\infty$. Let $|\bold{a}_i|=t_i$, in which $\bold{a}_i$ is the $b$-ary expansion of $m_i$. Our first result gives some comparison between the numbers $m_i$ and $b^{d_i}$.

\begin{lemma}\label{lem5}
There exist some constants $0<\underline{c}_1\leq \overline{c}_1<1$ , such that
$$b^{\underline{c}_1d_i}\leq m_i\leq  b^{\overline{c}_1d_i}$$
for any $1\leq i<\infty$. The two constants are independent of the script $i$.
\end{lemma} 
\begin{proof}
Since 
$$m_i=s_1^{r_{i,1}}s_2^{r_{i,2}}\cdots s_l^{r_{i,l}}\geq s_1^{d_i}=b^{d_i\log_{b}{s_1}}$$
and
$$m_i\leq s_l^{d_i}=b^{d_i\log_{b}{s_l}},$$
by simply taking $\underline{c}_1=\log_{b}{s_1}$ and  $\overline{c}_1=\log_{b}{s_l}$ we obtain the result.
\end{proof}

The next result concerns the comparison of the two numbers $t_i$ and $d_i$.
\begin{lemma}\label{lem6}
There exist some constants $0<\underline{c}_2\leq \overline{c}_2\leq 2$ , such that
$$\underline{c}_2d_i\leq t_i\leq  \overline{c}_2d_i$$
for any $1\leq i<\infty$. The two constants are independent of the script $i$.
\end{lemma} 
\begin{proof}
Note that $t_i=\lfloor\log_{b}{m_i}\rfloor+1$, so
$$t_i> \log_{b}{m_i}\geq \log_{b}{b^{\underline{c}_1d_i}}=\underline{c}_1d_i$$
while
$$t_i\leq \log_{b}{m_i}+1\leq \log_{b}{b^{\overline{c}_1d_i}}+1=\overline{c}_1d_i+1$$
for any $1\leq i<\infty$, in virtue of Lemma \ref{lem5}. Considering (\ref{eq27}) we have $\lim_{i\rightarrow\infty} d_i=\infty$. Then by taking $\underline{c}_2=\underline{c}_1$ and  $\overline{c}_2=\min\{\sup_i \{\overline{c}_1+\frac{1}{d_i}\},2\}$ we justify the result.
\end{proof}

\begin{rem}
In fact we can always choose  $\overline{c}_2<1$ in Lemma \ref{lem6}, except the case that the base $b$ is a prime and $m_1=b$.
\end{rem}

In the following we consider the overlapping component of the two numbers $b^{t_i}$ and $m_i$ for $i\in\mathbb{N}_+$.
\begin{lemma}\label{lem7}
There exists some constant $0<\underline{c}_3<1$ , such that
$$gcd(b^{t_i},m_i)\geq  b^{\underline{c}_3d_i}$$
for any $1\leq i<\infty$. The constant is independent of the script $i$.
\end{lemma} 
\begin{proof}
Since $m_i$ admits exactly $l$ prime factors, by the pigeonhole principle, there exists some $1\leq j\leq l$, such that
$$r_{i,j}\geq \frac{1}{l}d_i$$
for any $i\in\mathbb{N}_+$. Though we do not write explicitly it is possible that $j$ depends on $i$. Note that $s_j^{k_jt_i}\geq s_j^{k_j\underline{c}_2d_i}$ is a factor of $b^{t_i}$ (the inequality is due to Lemma \ref{lem6}). Let $\underline{c}_{3*}=\min\{\frac{1}{l}, k_j\underline{c}_2\}$. Then we have
$$gcd(b^{t_i},m_i)\geq s_j^{\underline{c}_{3*} d_i}=b^{(\log_{b}{s_j})\underline{c}_{3*} d_i}\geq b^{(\log_{b}{s_1})\underline{c}_{3*} d_i}$$
for any $1\leq i<\infty$.  Finally, by taking $\underline{c}_3=(\log_{b}{s_1})\underline{c}_{3*}$ we justify the result.
\end{proof}
  
Let $n_i=\sum_{j=1}^i t_j$. In the following we turn to properties of the numerators $\{u_{n_i}\}_{i\in\mathbb{N}_+}$ of the rational truncations of the number $(.\bold{a}_1\bold{a}_2\bold{a}_3\cdots)_b$. Considering the $S$-components of them, Lemma \ref{lem7} induces the following result instantly.

\begin{coro}\label{cor6}
The $S$-component of $u_{n_i}$ satisfies
$$c_{u_{n_i},S}\geq b^{\underline{c}_3d_i}$$
for any $1\leq i<\infty$, in which the constant $\underline{c}_3$ is from  Lemma \ref{lem7}.
\end{coro}   
\begin{proof}
In case $i=1$ the result is obvious. We mention that the following inductive formula 
$$u_{n_i}=b^{t_i} u_{n_{i-1}}+m_i$$
holds for any $2\leq i<\infty$. Since $F_b\cup F_{m_i}\subset S$,  we have $gcd(b^{t_i},m_i)|c_{u_{n_i},S}$. Then the result follows from Lemma \ref{lem7}. 

\end{proof}  

We also compare two numbers $u_{n_i}$ and $b^{d_i}$ in the following result, which may be of interests for other purposes.
\begin{pro}\label{lem8}
There exist some constants $0<\underline{c}_4\leq \overline{c}_4<\infty$ , such that
$$b^{\underline{c}_4d_i}\leq u_{n_i}\leq  b^{\overline{c}_4d_i}$$
for any $1\leq i<\infty$. The two constants are independent of the script $i$.
\end{pro}
\begin{proof}
The lower bound is obvious. Since $u_{n_i}\geq m_i$, in virtue of Lemma \ref{lem5}, we can simply take $\underline{c}_4=\underline{c}_1$. To build the upper bound, note that
\begin{equation}\label{eq29}
u_{n_i}=b^{t_i} u_{n_{i-1}}+m_i=b^{t_i} \sum_{j=1}^{i-1} b^{\sum_{q=j+1}^{i-1}t_q}m_j+m_i
\end{equation}
for any $1\leq i<\infty$ (certain summary terms in the formula disappear in case their upper scripts are less than the lower ones). Note that $b^{t_i}$ and 
$m_i$ are both well bounded from the above according to Lemma \ref{lem5} and \ref{lem6}, so we focus on the upper bound of the term $\sum_{j=1}^{i-1} b^{\sum_{q=j+1}^{i-1}t_q}m_j$ in the following. In this case we have

\begin{equation}\label{eq28}
\begin{array}{ll}
& \sum_{q=j+1}^{i-1}t_q\vspace{2mm}\\
\leq & \overline{c}_2\sum_{q=j+1}^{i-1}d_q\vspace{2mm}\\
\leq & \overline{c}_2\Big(\frac{1}{(1+\epsilon_{1.2})^{i-j-1}}d_i+\frac{1}{(1+\epsilon_{1.2})^{i-j-2}}d_i+\cdots+\frac{1}{(1+\epsilon_{1.2})^2}d_i+\frac{1}{1+\epsilon_{1.2}}d_i\Big)\vspace{2mm}\\
<&\cfrac{\overline{c}_2}{\epsilon_{1.2}}d_i\vspace{2mm}\\
\end{array}
\end{equation}
for any $1\leq j\leq i-2$ (in case $j=i-1$ the summation is empty). The first inequality is due to Lemma \ref{lem6}, the second one is due to (\ref{eq27}). Then
\begin{equation}\label{eq30}
\begin{array}{ll}
& \sum_{j=1}^{i-1} b^{\sum_{q=j+1}^{i-1}t_q}m_j\vspace{2mm}\\
\leq & \sum_{j=1}^{i-1} b^{\frac{\overline{c}_2}{\epsilon_{1.2}}d_i}m_j\vspace{2mm}\\
\leq & b^{\frac{\overline{c}_2}{\epsilon_{1.2}}d_i}\sum_{j=1}^{i-1} b^{\overline{c}_1d_j}\vspace{2mm}\\
\leq & b^{\frac{\overline{c}_2}{\epsilon_{1.2}}d_i} b^{\sum_{j=1}^{i-1}\overline{c}_1d_j}\vspace{2mm}\\
\leq & b^{\frac{\overline{c}_2}{\epsilon_{1.2}}d_i} b^{\overline{c}_1(\frac{1}{(1+\epsilon_{1.2})^{i-1}}d_i+\frac{1}{(1+\epsilon_{1.2})^{i-2}}d_i+\cdots+\frac{1}{(1+\epsilon_{1.2})^2}d_i+\frac{1}{1+\epsilon_{1.2}}d_i)}\vspace{2mm}\\
< & b^{(\frac{\overline{c}_2+\overline{c}_1}{\epsilon_{1.2}})d_i}. \vspace{2mm}\\
\end{array}
\end{equation}
The first inequality is due to (\ref{eq28}), the second one is due to Lemma \ref{lem5}, the fourth one is due to (\ref{eq27}). Now combing all the previous estimations Lemma \ref{lem5}, \ref{lem6} and (\ref{eq30}) through (\ref{eq29}), we have
$$u_{n_i}<b^{\overline{c}_2d_i}b^{(\frac{\overline{c}_2+\overline{c}_1}{\epsilon_{1.2}})d_i}+b^{\overline{c}_1d_i}<b^{(\overline{c}_2+\overline{c}_1)(1+\frac{1}{\epsilon_{1.2}})d_i}$$
for any $1\leq i<\infty$. Finally, taking $\overline{c}_4=(\overline{c}_2+\overline{c}_1)(1+\frac{1}{\epsilon_{1.2}})$ we arrive at the result.

\end{proof}
Equipped with all the estimations above, we are now in a position to justify Theorem \ref{thm5}.
\begin{proof}[Proof of Theorem \ref{thm5}]
For the number $(.\bold{a}_1\bold{a}_2\bold{a}_3\cdots)_b$, consider its rational truncations at $n_i=\sum_{j=1}^i t_j$ for $i\in\mathbb{N}_+$. Note that 
$$n_i=\sum_{j=1}^i t_j\leq \overline{c}_2\sum_{j=1}^i d_j<\overline{c}_2 (1+\frac{1}{\epsilon_{1.2}})d_i$$
for any $1\leq i<\infty$, so
\begin{equation}\label{eq31}
d_i\geq \frac{\epsilon_{1.2}}{\overline{c}_2(\epsilon_{1.2}+1)}n_i
\end{equation}
for any $1\leq i<\infty$. In virtue of (\ref{eq31}) and Corollary \ref{cor6}, the numerator of the $n_i$-th truncation and its $S$-component satisfy

$$\cfrac{c_{u_{n_i},S}}{u_{n_i}}\geq \frac{b^{\underline{c}_3d_i}}{b^{n_i}}\geq \frac{b^{\underline{c}_3\frac{\epsilon_{1.2}}{\overline{c}_2(\epsilon_{1.2}+1)}n_i}}{b^{n_i}}=b^{(\frac{\underline{c}_3\epsilon_{1.2}}{\overline{c}_2(\epsilon_{1.2}+1)}-1)n_i}$$
for any $1\leq i<\infty$. Finally the result follows from an application of Theorem \ref{thm1} with $\epsilon_{1.1}=\frac{\underline{c}_3\epsilon_{1.2}}{\overline{c}_2(\epsilon_{1.2}+1)}$.

\end{proof} 

In fact there are some flexibilities in Theorem \ref{thm5}, which can be manipulated to build more  transcendentals at certain demands.  

\begin{coro}
The requirement (\ref{eq27}) can be relaxed to
\begin{equation}\label{eq32}
\lim_{j\rightarrow\infty}\inf_{j\leq i<\infty}\Big\{\cfrac{\sum_{j=1}^l r_{i+1,j}}{\sum_{j=1}^l r_{i,j}}\Big\}>1
\end{equation}
while Theorem \ref{thm5} still holds. 
\end{coro}   
\begin{proof}
Assume (\ref{eq32}) holds, we can find some $j_*\in\mathbb{N}_+$ and some $\epsilon_{1.2}>0$ such that 
$$\inf_{j_*\leq i<\infty}\Big\{\cfrac{\sum_{j=1}^l r_{i+1,j}}{\sum_{j=1}^l r_{i,j}}\Big\}\geq 1+\epsilon_{1.2}.$$
An application of Theorem \ref{thm5} to the nodes $\{m_i=s_1^{r_{i,1}}s_2^{r_{i,2}}\cdots s_l^{r_{i,l}}\}_{i=j_*}^\infty$
guarantees that the number 
\begin{center}
$.\bold{a}_{j_*}\bold{a}_{j_*+1}\bold{a}_{j_*+2}\cdots$
\end{center}
is transcendental, which implies $.\bold{a}_1\bold{a}_2\bold{a}_3\cdots$ is transcendental too.
\end{proof} 

\begin{rem}
The condition (\ref{eq32}) in fact imposes at least exponential growth of the sequence of indexes $\{d_i\}_{i\in\mathbb{N}_+}$, while we suspect some slower growth rate will be enough to guarantee some similar conclusion as Theorem \ref{thm5}.  
\end{rem}

One can add words of relatively slower growing length between the words $\{\bold{a}_i\}_{i\in\mathbb{N}_+}$ in Theorem \ref{thm5} while the resulted number is still  transcendental.

\begin{theorem}\label{thm6}
Let $b\geq 2$ be some integer base with its prime factorization $b=s_1^{k_1}s_2^{k_2}\cdots s_l^{k_l}$ for some finite set of primes $S=\{s_1, s_2, \cdots, s_l\}$ with positive integers $\{k_1,k_2,\cdots,k_l,l\}$. Let 
$\{m_i=s_1^{r_{i,1}}s_2^{r_{i,2}}\cdots s_l^{r_{i,l}}\}_{i=1}^\infty$
be a sequence of nodes satisfying (\ref{eq27}) for some fixed $\epsilon_{1.2}>0$ and $d_i=\sum_{j=1}^l r_{i,j}$.  Then the number with $b$-ary expansion
\begin{center}
$\bold{a}_1\bold{a}_1'\bold{a}_2\bold{a}_2'\bold{a}_3\bold{a}_3'\cdots$
\end{center}
in $(0,1)$ is transcendental, in which $\bold{a}_i$ is the $b$-ary expansion of the node $m_i$ and $\{\bold{a}_i'\}_{i\in\mathbb{N}_+}$ is an arbitrary sequence of words satisfying
\begin{equation}\label{eq33}
\sum_{j=1}^{i-1} |\bold{a}_j'|=o(d_i)
\end{equation}
as $i\rightarrow\infty$. 
\end{theorem}
\begin{proof}
The proof follows in line with the proof of Theorem \ref{thm5}. The condition (\ref{eq33}) guarantees $\bold{a}_i$ ($m_i$) dominates the contribution of $\bold{a}_i'$ ($m_i'$ whose $b$-ary expansion being $\bold{a}_i'$), while all the estimations still hold for the extended word. Details are left to the enthusiastic readers.

\end{proof}

We build some new transcendentals in virtue of Theorem \ref{thm6} in the following.

\begin{exa}
Let $b=12=2^2 3$ and $\{m_i=3^{2^i}\}_{i\in\mathbb{N}}$ be a sequence of nodes satisfying (\ref{eq27}) with $\epsilon_{1.2}=1$.  Then according to Theorem \ref{thm5} the number
$$(.9\ 6\ 9\ 3\ 9\ 6\ 9\ 1\ 2\ 4\ 11\ 11\ 3\ 6\ 9\ 1\ 5\ 3\ 9\ 11\ 3\ 8\ 8\ 7\ 11\ 11\ 2\ 3\ 6\ 9\ ...)_{12}$$
is transcendental. Moreover, by adding 10 after every blocks successively in the above expansion, we are sure that the number
$$(.9\ \textcolor{red}{10}\ 6\ 9\ \textcolor{red}{10}\  3\ 9\ 6\ 9\ \textcolor{red}{10}\ 1\ 2\ 4\ 11\ 11\ 3\ 6\ 9\ \textcolor{red}{10}\ 1\ 5\ 3\ 9\ 11\ 3\ 8\ 8\ 7\ ...)_{12}$$
is still transcendental, in virtue of Theorem \ref{thm6}.
\end{exa}

We can also do some perturbations on the nodes in Theorem \ref{thm5}, while the resulted number preserves its  transcendence.

\begin{theorem}\label{thm7}
Let $b\geq 2$ be some integer base with its prime factorization $b=s_1^{k_1}s_2^{k_2}\cdots s_l^{k_l}$ for some finite set of primes $S=\{s_1, s_2, \cdots, s_l\}$ with positive integers $\{k_1,k_2,\cdots,k_l,l\}$. Let 
$\{m_i= s_1^{r_{i,1}}s_2^{r_{i,2}}\cdots s_l^{r_{i,l}}\}_{i=1}^\infty$
be a sequence of nodes satisfying (\ref{eq27}) for some fixed $\epsilon_{1.2}>0$ and $d_i=\sum_{j=1}^l r_{i,j}$.  Let $\{w_i\}_{i\in\mathbb{N}}$ be a sequence of positive integers satisfying  
\begin{equation}\label{eq34}
\log w_i=o(d_i)
\end{equation}
as $i\rightarrow\infty$. Then the number with $b$-ary expansion
\begin{center}
$\bold{a}_1\bold{a}_2\bold{a}_3\cdots$
\end{center}
in $(0,1)$ is transcendental, in which $\bold{a}_i$ is the $b$-ary expansion of the node $w_im_i$.
\end{theorem}

\begin{proof}
The condition (\ref{eq34}) guarantees the loss of ratio of the $S$-component of the numerators of the rational truncations of the concerned number over the numerators is well under control. Details are left to the enthusiastic readers.

\end{proof}

We can build new transcendentals in virtue of Theorem \ref{thm7} now.

\begin{exa}
Let $b=30=2\cdot3\cdot5$, $\{m_i=2^{2^i}3^i\}_{i\in\mathbb{N}}$ and $\{w_i=i\}_{i\in\mathbb{N}}$.  Then according to Theorem \ref{thm7} the number
$$(.12\ 9\ 18\ 23\ 1\ 6\ 26\ 6\ 12\ 28\ 24\ 7\ 28\ 18\ 8\ 10\ 28\ 9\ 18\ 0\ ...)_{30}$$
is transcendental.
\end{exa}

\section{Eventually periodic words are not finite factorial}\label{sec7}

Although we focus on the (partially) finite factorial property of irrationals (or non-eventually periodic words), we do wonder whether these properties apply to rationals (or eventually periodic words) essentially. It turns out that the finite factorial property is a nice indicator to distinguish  irrationals (none-eventually periodic words) from rationals (eventually periodic words).

We mean to prove Theorem \ref{pro3} in this section. We will achieve it through several steps.  The notion of \emph{coprime sequences} will be used in the following.

\begin{defn}
A sequence of (finite or infinite) integers $\{r_i\}_{i=1}^k$ of length $k\in\mathbb{N}_+\cup\{\infty\}$ is called a coprime sequence if $gcd(r_1,r_2,r_3,\cdots)=1$. It is called a coprime sequence up to an integer $r$ if $\{r_i/r\}_{i=1}^k$ is a coprime sequence.
\end{defn} 

To test the finite factorial property of a word (or a number) with respect to some base, we introduce some quantitive indicators. 
\begin{defn}
For some base $b\geq 2$ and some infinite word $\bold{a}\in \bold{A}_b$, the factorial index of the word $\bold{a}$ and a strictly increasing subsequence $\{n_i\}_{i=1}^\infty$ of $\mathbb{N}_+$ is defined as 
$$v_{\{n_i\},\bold{a}}=\lim_{j\rightarrow\infty}\#\cup_{i=1}^j F_{u_{n_i,\bold{a}}}.$$
 The factorial index of the word $\bold{a}$ is defined as 
$$v_{\bold{a}}=\liminf_{\{n_i\}_{i=1}^\infty\subset\mathbb{N}_+} v_{\{n_i\},\bold{a}}.$$
\end{defn} 
We will simply write it as $v_{\{n_i\}}$ in case the word is obvious.  It is easy to see that an infinite word is finite factorial if and only if its factorial index is finite (for certain subsequence $\{n_i\}_{i=1}^\infty$). Be careful that it is possible 
$$v_{\mathbb{N}_+,\bold{a}}=\infty$$
for some finite factorial word $\bold{a}\in B^\infty$. 

The following result indicates that the factorial indexes of words and subsequences are sensitive to initial letters in the words.
\begin{pro}
For some base $b\geq 2$, suppose the factorial index of some infinite word $\bold{a}=a_1a_2a_3\cdots\in \bold{A}_b$ and some strictly increasing subsequence $\{n_i\}_{i=1}^\infty$ of $\mathbb{N}_+$ satisfies
\begin{equation}\label{eq20}
v_{\{n_i\},\bold{a}}<\infty.
\end{equation}  
Now if its prefix $\bold{a}|_k\neq \bold{0}$ for some $k\in\mathbb{N}_+$, then
\begin{equation}\label{eq21}
v_{\{n_i\},\bold{0}\bold{a}_{-k}}=v_{\{n_i-k>0\},\bold{a}_{-k}}=\infty,
\end{equation}
in which $\bold{a}_{-k}=a_{k+1}a_{k+2}a_{k+3}\cdots$.
\end{pro} 
\begin{proof}
It is easy to see that 
\begin{equation}
u_{n_i,\bold{0}\bold{a}_{-k}}+b^{n_i-k}u_{k,\bold{a}}=u_{n_i,\bold{a}}
\end{equation}
for any $i$ such that $n_i>k$. Without loss of generality we can assume $u_{n_i,\bold{a}} \perp b$ for any $i$ such that $n_i>k$, then $u_{n_i,\bold{a}}$ is coprime with $u_{n_i,\bold{0}\bold{a}_{-k}}$ (posssibly up to a factor of $u_{k,\bold{a}}$). Note that $u_{k,\bold{a}}>0$ since $\bold{a}|_k\neq \bold{0}$. Now suppose (\ref{eq21}) does not hold, that is, 
$$v_{\{n_i\},\bold{0}\bold{a}_{-k}}<\infty.$$ 
Considering (\ref{eq20}), let
$$S=\big(\cup_{n_i>k} F_{u_{n_i,\bold{0}\bold{a}_{-k}}}\big)\cup\big(\cup_{n_i>k} F_{u_{n_i,\bold{a}}}\big)\cup F_b$$
be the finite set of primes, then every coprime sequence $(u_{n_i,\bold{0}\bold{a}_{-k}},b^{n_i-k},u_{n_i,\bold{a}})$ generated by $S$ gives a solution to the equation 
\begin{equation}\label{eq22}
x+yu_{k,\bold{a}}=z
\end{equation}  
for any $i$ such that $n_i>k$. This causes contradiction to the upper bound on the number of solutions of $(\ref{eq22})$ by Lewis-Mahler (\cite{LM}).

\end{proof}

We start from exploring  finite factorial property of the simplest periodic word in the following.
\begin{lemma}\label{pro4}
For any base $b\geq 2$, the infinite periodic word $\bold{1}=111\cdots$ is not finite factorial with respect to $b$.  
\end{lemma} 
\begin{proof}
We will show that for any subsequence $\{n_i\}_{i=1}^\infty$ of $\mathbb{N}_+$, its factorial index satisfies
\begin{equation}\label{eq12}
\lim_{i\rightarrow\infty} v_{\{n_i\}}=\infty,
\end{equation}
which is enough to induce the conclusion. We justify (\ref{eq12}) according to varying factorial properties of the sequences of subscripts $\{n_i\}_{i=1}^\infty$.

\begin{enumerate}
\item[Case I]: $\#\cup_{i=1}^\infty F_{n_i}=\infty$.

 In this case we can find an infinite sequence of mutually distinct primes $\{s_j\}_{j=1}^\infty\subset P$, such that 
\begin{equation}\label{eq13}
s_j|n_{i_j}
\end{equation} 
for any $j\in\mathbb{N}_+$, in which $\{n_{i_j}\}_{j\in\mathbb{N}_+}$ is a subsequence of $\{n_i\}_{i=1}^\infty$. Note that any two members in $\{u_{s_j}\}_{j=1}^\infty$ are coprime to each other. This implies (\ref{eq12}) since 
$$u_{s_j}|u_{n_{i_j}}$$  
in virtue of (\ref{eq13}).

\item[Case II]: $\#\cup_{i=1}^\infty F_{n_i}<\infty$.

This means that there exists some finite set $S$ of primes  such that $\{n_i\}_{i=1}^\infty\subset M_S$. Then we can find a strictly increasing sequence $\{r_j\}_{j=1}^\infty$ such that   
\begin{equation}\label{eq14}
s^{r_j}|n_{i_j}
\end{equation}
for some $s\in S$ and some subsequence $\{n_{i_j}\}_{j\in\mathbb{N}_+}$ of $\{n_i\}_{i=1}^\infty$. Let 
$$q_k=1+b^{s^{k-1}}+b^{2s^{k-1}}+b^{3s^{k-1}}+\cdots+b^{(s-1)s^{k-1}}$$
for $1\leq k<\infty$. We claim that $q_k$ is coprime with $\prod_{j=1}^{k-1} q_j$, possibly up to $s$, for any $k\geq 2$. To see this, first note that 
$$\prod_{j=1}^{k-1} q_j=1+b+\cdots+b^{(s-1)(1+s+\cdots+s^{k-2})}=1+b+\cdots+b^{s^{k-1}-1}=\cfrac{b^{s^{k-1}}-1}{b-1}$$   
for any $k\geq 2$. Then the claim follows from the equation
$$q_k=sb^{(s-1)s^{k-1}}-(b^{s^{k-1}}-1)\sum_{j=1}^{s-1} (s-j)b^{(s-j-1)s^{k-1}}$$ 
for any $k\geq 2$. The importance of the sequence $\{q_k\}_{k\in\mathbb{N}_+}$ is that
\begin{equation}\label{eq15}
\prod_{j=1}^k q_j=u_{s^k}
\end{equation}
for any $k\geq 2$. It turns out that our claim guarantees that the sequence of cardinalities $\{\# \cup_{1\leq i\leq j}F_{u_{s^{r_i}}}\}_{1\leq j<\infty}$ is strictly increasing with respect to $j\in\mathbb{N}_+$, considering (\ref{eq15}). Finally, since  
$$u_{s^{r_j}}|u_{n_{i_j}}$$  
in virtue of (\ref{eq14}), we get (\ref{eq12}).
\end{enumerate}
\end{proof}

\begin{rem}
In case II in the proof of Lemma \ref{pro4}, if the base $b=2$ and $s=2$, we have
$$q_k=1+2^{2^{k-1}}$$ 
for $k\in\mathbb{N}_+$. These numbers are called Fermat numbers (refer to \cite{BC} by Boklan and Conway). It is well known  that any two members in $\{1+2^{2^{k-1}}\}_{k\in\mathbb{N}_+}$ are coprime. Members in $\{1+2^{2^{k-1}}\}_{k=0}^4$ are all primes, while it is conjectured that the ones in $\{1+2^{2^{k-1}}\}_{k=5}^\infty$ are all composites. These numbers play important role in a problem of Landau (see for example \cite{Iwa} by Iwaniec). Our sequence $\{q_k\}$ is some generalised Fermat numbers in certain sense, whose primality may be of individual interest. 
\end{rem}

The following result is instant in virtue of Lemma \ref{pro4}.

\begin{coro}\label{cor3}
For any $1\leq k\leq b-1$, the infinite periodic word $\bold{k}=kkk\cdots$ is not finite factorial with respect to any base $b\geq 2$.
\end{coro}
\begin{proof}
The result follows from the fact $u_{i,\bold{k}}=ku_{i,\bold{1}}$ for any $i\in\mathbb{N}_+$ and Lemma \ref{pro4}.
\end{proof}

While elementary methods may still be used to decide the (none) finite factorial properties of more words in our expectation, we switch to Evertse's Theorem (\cite{Eve}) on the upper bounds of numbers of solutions of equations in $S$-units on arbitrarily algebraic number field to continue our explorations in the following.\footnote{We learn Evertse's Theorem from the work of Gy\"ory-Stewart-Tijdeman \cite{GST}, where they employ it to tackle a problem of Erd\"os and Tur\'an (\cite{ET}).}    

\begin{lemma}\label{lem4}
Any periodic word in $\bold{A}_{b}$ is not finite factorial with respect to any base $b\geq 2$.
\end{lemma}
\begin{proof}
Let $\bold{a}=a_1a_2a_3\cdots\in \bold{A}_{b}$ be a periodic word of essential period $r\in\mathbb{N}_+$. In case $r=1$ the result follows from Corollary \ref{cor3}, so in the following we assume $r\geq 2$. We justify the result by reduction to absurdity. Assume $\bold{a}$ is finite factorial, then there exist some strictly increasing subsequence $\{n_i\}_{i=1}^\infty$ of $\mathbb{N}_+$ and some finite set $S$ of primes, such that $u_{n_i}\in M_S$ for any $1\leq i<\infty$. By the Pigeon-hole Principle, we can find an infinite subsequence $\{n_{i_j}\}_{j=1}^\infty$, such that $n_{i_j}=q_j r+k$ for some fixed $0\leq k\leq r-1$ and some strictly increasing sequence $\{q_j\}_{j=1}^\infty\subset\mathbb{N}_+$. For the sequence of subscripts $\{n_{i_j}\}_{j=1}^\infty$, we have
$$u_{n_{i_j}}=u_rb^k(1+b^r+b^{2r}+\cdots+b^{(q_j-1)r})+u_k=\cfrac{u_rb^k(b^{q_jr}-1)}{b^r-1}+u_k$$  
for any $j\geq 1$. Then we get
\begin{equation}
(b^r-1)u_{n_{i_j}}=u_rb^{q_jr+k}+u_k(b^r-1)-u_rb^k
\end{equation}
for any $j\geq 1$. Note that $u_k(b^r-1)-u_rb^k>0$ for any $0\leq k\leq r-1$ since we are assuming $r\geq 2$. This means that every coprime triple $(b^{q_jr+k}, 1, u_{n_{i_j}})$ gives a solution to the equation 
\begin{equation}\label{eq16}
u_r x+\big(u_k(b^r-1)-u_r b^k\big)y=(b^r-1)z
\end{equation} 
for any $j\geq 1$. Now we consider the finite set of primes  $S\cup F_{b}$, note that $\{b^{q_jr}, u_{n_{i_j}}\}\subset M_{S\cup F_{b}}$. Applying the Evertse's Theorem we see that there can only be finitely many solutions generated by the set $S\cup F_{b}$, which contradicts the fact that (\ref{eq16}) holds with  $(b^{q_jr+k}, 1, u_{n_{i_j}})$ for any $j\geq 1$. 
\end{proof}

The case $k=0$ in the proof of Lemma \ref{lem4}  follows from Lemma \ref{pro4} directly in fact. The arguments apply to eventually periodic words essentially.

\begin{proof}[Proof of Theorem \ref{pro3}]
For some $r'\in\mathbb{N}_+$, let
$$\bold{a}=\bold{a}'\bold{a}''\cdots\in \bold{A}_{b}$$
be the concatenation of some finite word $\bold{a}'=a'_1a'_2a'_3\cdots a'_{r'}$ of length $r'\in\mathbb{N}_+$ and some infinite periodic word $\bold{a}''=a''_1a''_2a''_3\cdots \in\bold{A}_{b}$ of essential period $r''\in\mathbb{N}_+$. In case $\bold{a}'=00\cdots0$ the result follows from Lemma \ref{lem4}, so we assume $\bold{a}'\neq 00\cdots0$ in the following. Furthermore we can assume $\bold{a}\neq 00\cdots0kkk\cdots$ for some fixed letter $1\leq k\leq b-1$ since this again falls into Corollary \ref{cor3}. Be careful that We will take the  numerators of the finite word $\bold{a}'$ and the infinite words $\bold{a}, \bold{a}''$ respectively in the following. For example, 
$$u_{r',\bold{a}}=u_{r',\bold{a}'}=b^{r'-1} a_1'+b^{r'-2} a_2'+\cdots+a_{r'}',$$
while
$$u_{r'',\bold{a}''}=b^{r''-1} a_1''+b^{r''-2} a_2''+\cdots+a_{r''}''.$$

We again resort to the reduction to absurdity. Now  assuming $\bold{a}$ is finite factorial, then there exist some subsequence $\{n_i\}_{i=1}^\infty$ of $\mathbb{N}_+$ and some finite set $S$ of primes, such that $u_{n_i,\bold{a}}\in M_S$ for any $1\leq i<\infty$. By the Pigeon-hole Principle, we can find an infinite subsequence $\{n_{i_j}\}_{j=1}^\infty$, such that $n_{i_j}=r'+q_j r''+k$ for some fixed $0\leq k\leq r''-1$ and some strictly increasing sequence $\{q_j\}_{j=1}^\infty\subset\mathbb{N}_+$. Taking the $n_{i_j}$-th numerator of the word $\bold{a}$, we have
\begin{equation}\label{eq17}
\begin{array}{ll}
& u_{n_{i_j},\bold{a}}\vspace{3mm} \\ 
=&b^{q_j r''+k}u_{r',\bold{a}'}+u_{r'',\bold{a}''}b^k(1+b^{r''}+b^{2r''}+\cdots+b^{(q_j-1)r''})+u_{k,\bold{a}''}\vspace{3mm}\\ 
=&b^{q_j r''+k}u_{r',\bold{a}'}+\cfrac{u_{r'',\bold{a}''}b^k(b^{q_jr''}-1)}{b^{r''}-1}+u_{k,\bold{a}''}
\end{array}
\end{equation} 
for any $j\geq 1$.   Rewriting (\ref{eq17}) as
$$
(b^{r''}-1)u_{n_{i_j},\bold{a}}=b^{q_jr''+k}(u_{r',\bold{a}'}(b^{r''}-1)+u_{r'',\bold{a}''})+u_{k,\bold{a}''}(b^{r''}-1)-b^ku_{r'',\bold{a}''}
$$
for any $j\geq 1$. Obviously we have $gcd(b^{q_jr''+k}, u_{n_{i_j},\bold{a}})< b^{r''}$ for any $j\geq 1$. Note that $u_{k,\bold{a}''}(b^{r''}-1)-b^ku_{r'',\bold{a}''}>0$ for any fixed $0\leq k\leq r''-1$. Then every coprime triple $(b^{q_jr''+k}, 1, u_{n_{i_j},\bold{a}})$ gives a solution to the equation 
\begin{equation}\label{eq18}
\big(u_{r',\bold{a}'}(b^{r''}-1)+u_{r'',\bold{a}''}\big) x+\big(u_{k,\bold{a}''}(b^{r''}-1)-b^ku_{r'',\bold{a}''}\big)y=(b^{r''}-1)z
\end{equation} 
for any $j\geq 1$. Considering $\{b^{q_jr''+k}, u_{n_{i_j},\bold{a}}\}\subset M_{S\cup F_{b}}$ and the Evertse's Theorem, we see that there can only be finitely many solutions generated by the set $S\cup F_{b}$ for the equation (\ref{eq18}), which contradicts the fact that (\ref{eq18}) holds with  $(b^{q_jr''+k}, 1, u_{n_{i_j},\bold{a}})$ for any $j\geq 1$.

\end{proof}

It would be an interesting project to explore the (none) partially finite factorial property of eventually periodic words in $\bold{A}_{b}$ with respect to some base $b\geq 2$.

\section{The set of finite factorial numbers}\label{sec3}

For a finite set $S=\{s_1, s_2, \cdots, s_l\}$ of $l$ primes and an integer base $b\geq 2$,  let 
\begin{center}
$
\begin{array}{ll}
& \Xi^{FF}_{S,b}\\
=&\{\xi=(.\bold{a})_b=(.a_1a_2a_3\cdots)_b: \bold{a} \mbox{\ is a finite factorial word with respect to the base\ } b\},\\
\end{array}
$
\end{center}
be the collection of all the FF numbers with respect to the $S, b$. Alternatively, let 
\begin{center}
$
\Xi^{PFF}_{S,b}=\{\xi=(.\bold{a})_b: \bold{a} \mbox{\ is a partially finite factorial word with respect to the base\ } b\},
$
\end{center}
be the collection of all the PFF numbers with respect to the $S, b$. Let 
$$\Xi^{FF}=\cup_{S\subset P, \# S<\infty, b\in\mathbb{N}_+\setminus\{1\}} \Xi^{FF}_{S,b}$$
and 
$$\Xi^{PFF}=\cup_{S\subset P, \# S<\infty, b\in\mathbb{N}_+\setminus\{1\}} \Xi^{PFF}_{S,b}$$
be the collections of all the finite factorial and partially finite factorial numbers respectively. Let 
\begin{center}
$M_S=\{m_{r_1, r_2, \cdots, r_l}=\prod_{k=1}^l s_k^{r_k}:  (r_1, r_2, \cdots, r_l)\in\mathbb{N}^l\}$
\end{center}
be the set of integers generated by the set $S$, which are termed \emph{nodes} in the following. 

Let $$\{f_a(n)=bn+a\}_{a=0}^{b-1}$$
be the family of  endomorphisms on $\mathbb{N}_+$ with respect to some base $b$.  Let 
$$g(n)=\lfloor\cfrac{n}{b}\rfloor$$
be their pull-back map, which satisfies $g\circ f_a=Id$ for any $0\leq a\leq b-1$. In the following we establish some relationship on $M_S$, although most concepts can be generalised onto $\mathbb{N}_+$ in fact.  For two nodes $m_1, m_2\in M_S$, we say that $m_1$ is the \emph{ancester} of $m_2$ with respect to some base $b\geq 2$ if there exists some finite word $a_1a_2\cdots a_k\in B^k$ of length $k\in\mathbb{N}$ such that
\begin{equation}\label{eq7}
f_{a_1a_2\cdots a_k}(m_1)=f_{a_k}\circ f_{a_{k-1}}\circ\cdots\circ f_{a_1}(m_1)=m_2,
\end{equation}
which is equivalent to the condition 
\begin{equation}\label{eq8}
g^k(m_2)=m_1.
\end{equation}
We write this relationship as $m_1\sim m_2$. It is not an equivalent relationship as $\sim$ does not obey the symmetric principle.  In case (\ref{eq7}) holds such that $f_{a_1a_2\cdots a_i}(m_1)\notin M_S$ for any $1\leq i\leq k-1$ we say $m_1$ is the \emph{parent} of $m_2$ (writing as $m_1\stackrel{par}{\sim} m_2$).  Alternatively one can say that $m_2$ is a \emph{descendant} (\emph{child}) of $m_1$ with respect to some base $b\geq 2$ if (\ref{eq7}) holds (with $ f_{a_1a_2\cdots a_i}(m_1)\notin M_S$ for any $1\leq i\leq k-1$).    A sequence of integers $\{m_i\}_{i=1}^k\subset M_S$ (or $\{m_i\}_{i=1}^\infty\subset M_S$) is called a \emph{chain} of length $k$ (or $\infty$) if it satisfies
\begin{center}
$m_i\stackrel{par}{\sim} m_{i+1}$,
\end{center}  
for any $1\leq i\leq k-1$ (or $1\leq i<\infty$).

\begin{defn}
For some base $b\geq 2$, a subset $N\subset\mathbb{N}_+$ is said to be \emph{extendable} in base $b$ if for any integer $m\in\mathbb{N}_+$, there exist some $i\in\mathbb{N}_+$ and $n\in N$, such that $b^i m\leq n\leq b^i m+b^i-1$.
\end{defn}

The notion is introduced by Bird (\cite{Bir}) in order to check the primality of large integers given in some base-$b$ system. He proved that $N$ is extendable in base $b$ if and only if the set $\Big\{(\cfrac{\log n}{\log b})\Big\}_{n\in N}$ is dense in the unit interval $(0,1)$. He also gave some sufficient conditions for a set to be extendable, one of which (refer to Bird's Theorem in the appendix) turns out to be a powerful tool for us to explore the  descendants (children) of nodes in $M_S$.\footnote{We run into the notion and Bird's work through Kuipers and Niederreiter's book \cite{KN} in the exercises therein. They are also related to the uniform distribution theory of sequences of numbers \emph{mod} $1$, which is a  topic dating back to Weyl (\cite{Wey}). In fact we conjecture that 
$$\lim_{n\rightarrow\infty}\cfrac{\#(M_S\cap\cup_{i=1}^n [b^im, b^im+b^i-1])}{\# (M_S\cap[1,b^{k+n}-1])}=\cfrac{1}{b^k}$$ for any node $m\in M_S$ whose $b$-ary expansion is of length $k$, in case $\#(S\cup F_b)>1$.}  

\begin{lemma}\label{pro5}
Let $S=\{s_1, s_2, \cdots, s_l\}$ be a finite set of primes and $b\geq 2$ be some base, such that $\#(S\cup F_b)>1$. Then any node in $M_S$ admits infinitely many descendants and children, while any node in $M_S$ admits at most finitely  many (possibly no) ancesters and at most one (possibly no) parent in $M_S$.
\end{lemma}
\begin{proof}
We only show the first half of the result, while the second half is obvious. For any node $m\in M_S$, let $\bold{a}_1=a_{1,1}a_{1,2}\cdots a_{1,k},\in B^k$ be its $b$-ary expansion, so we have 
$$u_{\bold{a}_1}=m.$$
We consider the descendants of $m$ in the following.   Since $\#(S\cup F_b)>1$, we can always choose some prime number $s\in S$ such that 
$$b\notin M_{\{s\}}=\{s^r\}_{r=0}^\infty\subset M_S.$$ 

Now for the sequence $\{s^r\}_{r=0}^\infty$,  in virtue of the Bird's Theorem, the first condition is satisfied, so $M_{\{s\}}$ is  extendable in base $b$.  We claim that 
\begin{equation}\label{eq19}
\#(M_{\{s\}}\cap\cup_{i=1}^\infty [b^im, b^im+b^i-1])=\infty
\end{equation}
since $\#(M_{\{s\}}\cap\cup_{i=1}^\infty [b^im, b^im+b^i-1])<\infty$ will prevent $M_{\{s\}}$ to be extendable in base $b$. (\ref{eq19}) implies that $m$ is necessary to admit infinitely many descendants. Finally, there must be infinitely many children among these descendants, since for any finite word $\bold{a}_2\in\bold{A}_b^F$ and any fixed $j\in \mathbb{N}_+$, the numerators
$$\{u_{\bold{a}_1\bold{a}_2\bold{a}_3}: \bold{a}_3\in B^j\}$$
can not all be descendants of $m$.

\end{proof}

We present the readers some explicit numerical results to illustrate Lemma \ref{pro5} in Example \ref{exa2} in Section \ref{sec5}. In case $\#(S\cup F_b)=1$, that is, $S=\{s_1\}$ while $b=s_1^k$ for some $k\in\mathbb{N}_+$, we have $M_S=\{s_1^r\}_{r=0}^\infty$. In this case any member in $M_S=\{s_1^r\}_{r=0}^\infty$ admits exactly one child and infinitely many descendants, while any member admits one parent except the ones in $\{s_1^r\}_{r=0}^{k-1}$. There are in fact totally $k$ different infinite chains in this case.  Lemma \ref{pro5} induces the uncountability of the set of FF numbers. 

\begin{coro}\label{cor2}
For any finite set $S=\{s_1, s_2, \cdots, s_l\}$ of $l$ primes and any integer base $b\geq 2$ satisfying $\#(S\cup F_b)>1$, the set $\Xi^{FF}_{S,b}$ is uncountably infinite. 
\end{coro}
\begin{proof}
Let $M_{S,0}$ be the collections of members in $M_S$ admitting no ancestors. Obviously we have $M_{S,0}\neq\emptyset$ (in fact one can show that $\#M_{S,0}=\infty$ under the condition $\#(S\cup F_b)>1$). Now in virtue of Lemma  \ref{pro5}, every member in $M_{S,0}$ admits infinitely many children, let 
$$M_{S,1}\subset M_S$$
be the collection of these children. Again in virtue of Lemma  \ref{pro5}, every member in $M_{S,1}$ admits infinitely many children, let 
$$M_{S,2}\subset M_S$$
be the collection of these children. Successively, every member in  $M_{S,i}$ admits infinitely many children in $M_{S,i+1}$, in virtue of Lemma  \ref{pro5}, for any $2\leq i<\infty$. So we get uncountably many infinite chains starting with members in  $M_{S,0}$, such that every chain results in some real number (whose $b$-ary expansion begins with some none-$\bold{0}$ letter) distinct from one another in $\Xi^{FF}_{S,b}$.    
\end{proof}

It is easy to see that $M_{S,i}\cap M_{S,j}=\emptyset$ for any $i\neq j\in\mathbb{N}$ and $M_S=\cup_{i=0}^\infty M_{S,i}$ in the proof of Corollary \ref{cor2}, that is, $\{M_{S,i}\}_{i=0}^\infty$ is a partition of $M_S$.  By adding reals which are multiples of the resulted ones with rationals in $\{\cfrac{1}{b^k}\}_{k=1}^\infty$, we in fact exhaust the finite factorial numbers in $\Xi^{FF}_{S,b}$ for fixed $S, b$ in the proof of Corollary \ref{cor2}. 

\begin{coro}\label{cor4}
The set $\Xi^{FF}$ is an uncountably infinite set.
\end{coro}
\begin{proof}
Note that 
\begin{equation}\label{eq23}
\Xi^{FF}=\cup_{i\in\mathbb{N}_+, b\in\mathbb{N}_+\setminus\{1\}} \Xi^{FF}_{P_i,b},
\end{equation}
while every set $\Xi^{FF}_{P_i,b}$ with $\#(P_i\cup F_b)>1$ is an  uncountably infinite set in virtue of Corollary \ref{cor2}. 
\end{proof}

The notion of descendant (child) in fact applies to finite words.

\begin{defn}
Let $S$ be a finite set of primes. For any finite word $\bold{a}_1=a_{1,1}a_{1,2}\cdots a_{1,|\bold{a}_1|}\in A_b^F$ with respect to some base $b\geq 2$, a node $m\in M_S$ is called a descendant of the word $\bold{a}_1$ if there exists some finite word $\bold{a}_2=a_{2,1}a_{2,2}\cdots a_{2,|\bold{a}_2|}$, such that $u_{\bold{a}_1\bold{a}_2}\in M_S$. A descendant is called a child of  the word $\bold{a}_1$ if there are no other descendants of the word which arise as the ancestors of $m$.  
\end{defn}

If $\#(S\cup F_b)> 1$, any finite word admits  infinitely many descendants and children, according to an analogous argument as the proof of Lemma \ref{pro5}.

\begin{pro}\label{thm2}
Let $S=\{s_1, s_2, \cdots, s_l\}$ be some finite set of primes and $b\geq 2$ be some base satisfying $\#(S\cup F_b)>1$ and $F_b\setminus S\neq \emptyset$. For any non-zero finite word $\bold{a}_0\in A_b^F$,  there exist some unique finite factorial number $\xi_{S,b,\bold{a}_0}^{FF,-}=.\bold{a}$ (with $\bold{a}\in\bold{A}_b$) and some strictly increasing subsequence $\{n_i>|\bold{a}_0|\}_{i=1}^\infty$ of $\mathbb{N}_+$, such that 
\begin{enumerate}[(1).]
\item $\bold{a}|_{|\bold{a}_0|}=\bold{a}_0$.

\item $F_{u_{n_i}}\subset S$ for any $1\leq i< \infty$.

\item $u_{n_1}$ is the least child of the word $\bold{a}_0$.

\item $u_{n_i}$ is the least child of the word $\bold{a}|_{n_{i-1}}$ for any $i\geq 2$.
\end{enumerate}
 
\end{pro}

\begin{proof}
The $b$-ary expansion of the finite factorial number with the required properties can be obtained by successively adding finite words starting from $\bold{a}_0$. First, since $\bold{a}_0$ admits infinitely many children, we choose the least one $m_1=s_1^{r_{1,1}}s_1^{r_{1,2}}\cdots s_l^{r_{1,l}}\in M_S$ for some $\bold{r}_1=(r_{1,1}, r_{1,2}, \cdots, r_{1,l})\in \mathbb{N}^l$, such that 
$$m_1=u_{\bold{a}_0\bold{a}_1}$$
for some finite word $\bold{a}_1\in A_b^F$. Now according to Lemma \ref{pro5}, $m_1$ admits infinitely many children, so we can choose the least one $m_2=s_1^{r_{2,1}}s_1^{r_{2,2}}\cdots s_l^{r_{2,l}}\in M_S$ for some $\bold{r}_2=(r_{2,1}, r_{2,2}, \cdots, r_{2,l})\in \mathbb{N}^l$, such that 
$$m_2=u_{\bold{a}_0\bold{a}_1\bold{a}_2}$$
for some finite word $\bold{a}_2\in A_b^F$. Then successively we can find the least child $m_i=s_1^{r_{i,1}}s_1^{r_{i,2}}\cdots s_l^{r_{i,l}}\in M_S$ for some $\bold{r}_i=(r_{i,1}, r_{i,2}, \cdots, r_{i,l})\in \mathbb{N}^l$ and some finite word $\bold{a}_i\in A_b^F$, such that 
$$m_i=u_{\bold{a}_0\bold{a}_1\bold{a}_2\cdots\bold{a}_i}$$
for any $i\geq 3$. Then  we obtain an infinite word 
$$\bold{a}=\bold{a}_0\bold{a}_1\bold{a}_2\cdots$$
which results in some finite factorial number $\xi_{S,b,\bold{a}_0}^{FF,-}=.\bold{a}$. Finally, let 
$$n_i=\sum_{j=0}^i|\bold{a}_j|$$
for any $i\geq 1$. The number $\xi_{\bold{a}}$ and the sequence $\{n_i\}_{i=1}^\infty$ satisfy all the four requirements of the proposition. 
\end{proof}

\begin{rem}
According to Corollary \ref{cor1} and Theorem \ref{pro3}, the resulted number $\xi_{S,b,\bold{a}_0}^{FF,-}$ must be transcendental for any $S,b$ with $\#(S\cup F_b)>1, F_b\setminus S\neq \emptyset$ and any finite word $\bold{a}_0$.
\end{rem}

For fixed finite word $\bold{a}_0$, $\xi_{S,b,\bold{a}_0}^{FF,-}$ is the least finite factorial number whose $b$-ary expansion begins with $\bold{a}_0$, in some sense. In case $F_b\setminus S= \emptyset$, that is, $F_b\subset S$,  for any non-zero finite word $\bold{a}_0\in A_b^F$, there is still some unique number $.\bold{a}$ and some strictly increasing subsequence $\{n_i>|\bold{a}_0|\}_{i=1}^\infty$ satisfying all the requirements of Proposition \ref{thm2}, however, the number will be a rational one while $\bold{a}\in\bold{A}_{b,0}$. One can refer to Example \ref{exa3} for some initial blocks of the expansion of the number $\xi_{S,b,\bold{a}_0}^{FF,-}$ with  respect to $S=\{2,5,7,11\}, b=3$ and $\bold{a}_0=1$.

In the following we gear towards computing the Hausdorff dimension of the set of all FF numbers, we start from dealing with the dimension of its subsets.

\begin{lemma}\label{lem2}
For any finite set $S=\{s_1, s_2, \cdots, s_l\}$ of $l$ primes and any integer base $b\geq 2$, we have
$$\dim_H \Xi^{FF}_{S,b}=0.$$ 
\end{lemma}
\begin{proof}
First note that the set $\Xi^{FF}_{S,b}$ admits the following decomposition,
\begin{equation}\label{eq9}
\Xi^{FF}_{S,b}=\cap_{i=1}^\infty\cup_{k=i}^\infty \Xi^{FF}_{S,b,k}=\cap_{i=1}^\infty\cup_{k=i}^\infty\{\xi=(.a_1a_2a_3\cdots)_b: u_{k,\xi}\in M_S\}\footnote{The inspiration comes from Good, cf. \cite[(2$\cdot$2)]{Goo}.}
\end{equation}
in which $\Xi^{FF}_{S,b,k}=\{\xi=(.a_1a_2a_3\cdots)_b: u_{k,\xi}\in M_S\}$. We claim that the cardinality 
$$\#\{\bold{a}_k=a_1a_2\cdots a_k\in B^k: u_{\bold{a}_k}\in M_S\}$$
grows at most polynomially with respect to $k\in\mathbb{N}_+$. To see this, note that 
$$\max\{u_{\bold{a}_k}: \bold{a}_k=a_1a_2\cdots a_k\in B^k\}=b^k-1$$ 
for any $k$. Recall that $s_1\geq 2$ is the least number in $S$. Then  
$$\#\{\bold{a}_k=a_1a_2\cdots a_k\in B^k: u_{\bold{a}_k}\in M_S\}\leq \big(k\cfrac{\log b}{\log s_1}\big)^l$$
for any $k$. It then follows that
$$H^\alpha(\Xi^{FF}_{S,b})=\lim_{i\rightarrow\infty} \sum_{k=i}^\infty \Big(\big(k\cfrac{\log b}{\log s_1}\big)^l \cfrac{1}{b^{\alpha k}}\Big)=0$$
for any $\alpha>0$. This is enough to force the Hausdorff dimension of the set $\Xi^{FF}_{S,b}$ to be null.

\end{proof}

Note that the set $\Xi^{FF}_{S,b,k}=\cup_{u_{a_1a_2\cdots a_k}\in M_S, a_1a_2\cdots a_k\in B^k}I_{a_1a_2\cdots a_k}$ is in fact a finite union of closed intervals. Moreover, since the terminals of the fundamental intervals are all rationals (of eventually periodic expansions), so (\ref{eq9}) can be written as
$$\Xi^{FF}_{S,b}=\cap_{i=1}^\infty\cup_{k=i}^\infty (\Xi^{FF}_{S,b,k})^o=\cap_{i=1}^\infty\cup_{k=i}^\infty \cup_{u_{a_1a_2\cdots a_k}\in M_S}I_{a_1a_2\cdots a_k}^o$$  
in fact, in virtue of Theorem \ref{pro3}. This implies that $\Xi^{FF}_{S,b}$ is a $G_\delta$ set.

Although $\Xi^{FF}_{S,b}$ for given set of primes $S$ and base $b\geq 2$ (or $\Xi^{FF}$) can be negligible from the viewpoint of dimension (of course measure) in the transcendental numbers, it carries amazing topological property, in the sense that every other transcendental number (in fact every real number) can be approximated by sequences of numbers in $\Xi^{FF}_{S,b}$ (or $\Xi^{FF}$).

\begin{lemma}\label{pro2}
For any finite set $S=\{s_1, s_2, \cdots, s_l\}$ of primes and any integer base $b\geq 2$ satisfying $\#(S\cup F_b)>1$, the set $\Xi^{FF}_{S,b}$ is a dense subset of the unit interval $[0,1]$. 
\end{lemma}
\begin{proof}
It suffices for us to show that for any open interval $(\xi_1,\xi_2)\subset [0,1]$ with $0\leq \xi_1<\xi_2\leq 1$, 
$$(\xi_1,\xi_2)\cap \Xi^{FF}_{S,b}\neq \emptyset$$
for given $S$ and $b$. To see this, we choose some finite word $\bold{a}_0=a_1a_2\cdots a_{i_0}$ of length $i_0\in\mathbb{N}_+$ such that $I_{\bold{a}_0}\subset (\xi_1,\xi_2)$ (such word exists for $i_0$ large enough). Then following similar steps as the proof of Lemma  \ref{thm2} we can repetitively add blocks of finite words $\bold{a}_k$ of length $i_k$ such that the resulted infinite word  $\bold{a}=\bold{a}_0 \bold{a}_1\bold{a}_2\cdots$ satisfies 
\begin{center}
$F_{u_{\bold{a}_0 \bold{a}_1\bold{a}_2\cdots\bold{a}_k}}\subset S$
\end{center}
for any $k\in\mathbb{N}_+$. This guarantees $\xi_{\bold{a}}=.\bold{a}\in\Xi^{FF}_{S,b}$.
\end{proof}

\begin{coro}\label{cor5}
The set $\Xi^{FF}$ is a dense subset of the unit interval $[0,1]$.
\end{coro}
\begin{proof}
This follows from the fact (\ref{eq23}) and Lemma \ref{pro2}.
\end{proof}

Now we are in a position to prove Theorem \ref{thm3}.

\begin{proof}[Proof of Theorem \ref{thm3}]
In virtue of Corollary \ref{cor4} and \ref{cor5}, we only need to show that $dim_H \Xi^{FF}=0$. Considering (\ref{eq23}), this follows from the countable stability of Hausdorff dimension and Lemma \ref{lem2} in fact.
\end{proof}

\section{Numerical exhibitions}\label{sec5}

We provide the readers some concrete examples to illustrate some of our results in previous sections, which are expected to aid ones to grasp ideas in corresponding results. The first one is on the ancestor-descendant (parent-child) relationship (which should be understood along with Lemma \ref{pro5}).

\begin{exa}\label{exa2}
Let $S=\{2\}$ and $b=6$, in this case $M_S=\{2^r\}_{r=0}^\infty$. We list the finite chains containing all the numbers $\{2^r\}_{r=0, 1, 2, \cdots, 54, 69, 72, 93, 98}\subset M_S$ in the following. Some intermediate points are not revealed in case they are not nodes for tidiness. 

Word: $52$,
\begin{center}
\begin{tikzcd}
\ \ \ 5 \arrow[r, "f_2"] & 2^5 \arrow[l, bend left, "g"].
\end{tikzcd}
\end{center}

Word: $144$,
\begin{center}
\begin{tikzcd}
1 \arrow[r,  "f_4"] & 10 \arrow[r,  "f_4"] & 2^6 \arrow[ll, bend left, "g^2"].
\end{tikzcd}
\end{center}

Word: $332$,
\begin{center}
\begin{tikzcd}
3 \arrow[r,  "f_3"] & 21 \arrow[r,  "f_2"] & 2^7 \arrow[ll, bend left, "g^2"].
\end{tikzcd}
\end{center}

Word: $1104$,
\begin{center}
\begin{tikzcd}
1 \arrow[r,  "f_1"] & 7 \arrow[r,  "f_0"] & 42 \arrow[r,  "f_4"] & 2^8.
\end{tikzcd}
\end{center}

Word: $2212$,
\begin{center}
\begin{tikzcd}
2^1 \arrow[r,  "f_2"] & 14 \arrow[r,  "f_1"] & 85 \arrow[r,  "f_2"] & 2^9.
\end{tikzcd}
\end{center}

Word: $4424$,
\begin{center}
\begin{tikzcd}
2^2 \arrow[r,  "f_4"] & 28 \arrow[r,  "f_2"] & 170 \arrow[r,  "f_4"] & 2^{10}.
\end{tikzcd}
\end{center}

Word: $13252$,
\begin{center}
\begin{tikzcd}
1 \arrow[r,  "f_3"] & 9 \arrow[r,  "f_2"] & 56 \arrow[r,  "f_5"] & 341 \arrow[r,  "f_2"] & 2^{11}.
\end{tikzcd}
\end{center}

Word: $30544$,
\begin{center}
\begin{tikzcd}
3 \arrow[r,  "f_0"] & 18 \arrow[r,  "f_5"] & 113 \arrow[r,  "f_4"] & 682 \arrow[r,  "f_4"] & 2^{12}.
\end{tikzcd}
\end{center}

Word: $101532$,
\begin{center}
\begin{tikzcd}
1 \arrow[r,  "f_0"] & 6 \arrow[r,  "f_1"] & 37 \arrow[r,  "f_5"] & 227 \arrow[r,  "f_3"] & 1365 \arrow[r,  "f_2"] & 2^{13}.
\end{tikzcd}
\end{center}

Word: $203504$,
\begin{center}
\begin{tikzcd}
2^1 \arrow[r,  "f_0"] & 12 \arrow[r,  "f_3"] & 75 \arrow[r,  "f_5"] & 455 \arrow[r,  "f_0"] & 2730 \arrow[r,  "f_4"] & 2^{14}.
\end{tikzcd}
\end{center}

Word: $411412$,
\begin{center}
\begin{tikzcd}
2^2 \arrow[r,  "f_1"] & 25 \arrow[r,  "f_1"] & 151 \arrow[r,  "f_4"] & 910 \arrow[r,  "f_1"] & 5461 \arrow[r,  "f_2"] & 2^{15}.
\end{tikzcd}
\end{center}

Word: $1223224$,
\begin{center}
\begin{tikzcd}
1 \arrow[r,  "f_2"] & 2^3 \arrow[r,  "f_2"] & 50 \arrow[r,  "f_3"] & 303 \arrow[r,  "f_2"] & 1820 \arrow[r,  "f_2"] & 10922 \arrow[r,  "f_4"] & 2^{16}.
\end{tikzcd}
\end{center}

Word: $2450452$,
\begin{center}
\begin{tikzcd}
2^1 \arrow[r,  "f_4"] & 2^4 \arrow[r,  "f_5"] & 101 \arrow[r,  "f_0"] & 606 \arrow[r,  "f_4"] & 3640 \arrow[r,  "f_5"] & 21845 \arrow[r,  "f_2"] & 2^{17}.
\end{tikzcd}
\end{center}

Word: $5341344$,
\begin{center}
\begin{tikzcd}
5 \arrow[r,  "f_3"] & 33 \arrow[r,  "f_4"] & 202 \arrow[r,  "f_1"] & 1213 \arrow[r,  "f_3"] & 7281 \arrow[r,  "f_4"] & 43690 \arrow[r,  "f_4"] & 2^{18}.
\end{tikzcd}
\end{center}

Word: $15123132$,
\begin{center}
\begin{tikzcd}
1 \arrow[r,  "f_5"] & 11 \arrow[r,  "f_1"] & 67 \arrow[r,  "f_2"] & 404 \arrow[r,  "f_3"] & 2427 \arrow[r,  "f_1"] & 14563 \arrow[r,  "f_3"] & 87381 \arrow[r,  "f_2"] & 2^{19}.
\end{tikzcd}
\end{center}

Word: $34250304$,
\begin{center}
\begin{tikzcd}
3 \arrow[r,  "f_4"] & 22 \arrow[r,  "f_2"] & 134 \arrow[r,  "f_5"] & 809 \arrow[r,  "f_0"] & 4854 \arrow[r,  "f_3"] & 29127 \arrow[r,  "f_0"] & 174762 \arrow[r,  "f_4"] & 2^{20}.
\end{tikzcd}
\end{center}

Word: $112541012$,
\begin{center}
\begin{tikzcd}[column sep=2cm, row sep=0.1cm]
1 \arrow[r,  "f_{12541012}"]  & 2^{21}.
\end{tikzcd}
\end{center}

Word: $225522024$,
\begin{center}
\begin{tikzcd}[column sep=2cm, row sep=0.1cm]
2^1 \arrow[r,  "f_{25522024}"] & 2^{22}.
\end{tikzcd}
\end{center}

Word: $455444052$,
\begin{center}
\begin{tikzcd}[column sep=2cm, row sep=0.1cm]
2^2 \arrow[r,  "f_{55444052}"] & 2^{23}.
\end{tikzcd}
\end{center}

Word: $1355332144$,
\begin{center}
\begin{tikzcd}[column sep=2cm, row sep=0.1cm]
1 \arrow[r,  "f_{355332144}"] & 2^{24}.
\end{tikzcd}
\end{center}

Word: $3155104332$,
\begin{center}
\begin{tikzcd}[column sep=2cm, row sep=0.1cm]
3 \arrow[r,  "f_{155104332}"] & 2^{25}.
\end{tikzcd}
\end{center}

Word: $10354213104$,
\begin{center}
\begin{tikzcd}[column sep=2cm, row sep=0.1cm]
1 \arrow[r,  "f_{0354213104}"] & 2^{26}.
\end{tikzcd}
\end{center}

Word: $21152430212$,
\begin{center}
\begin{tikzcd}[column sep=2cm, row sep=0.1cm]
2^1 \arrow[r,  "f_{1152430212}"] & 2^{27}.
\end{tikzcd}
\end{center}

Word: $42345300424$,
\begin{center}
\begin{tikzcd}[column sep=2cm, row sep=0.1cm]
2^2 \arrow[r,  "f_{2345300424}"] & 2^{28}.
\end{tikzcd}
\end{center}

Word: $125135001252$,
\begin{center}
\begin{tikzcd}[column sep=2cm, row sep=0.1cm]
1 \arrow[r,  "f_2"] & 2^3 \arrow[r,  "f_{5135001252}"] & 2^{29}.
\end{tikzcd}
\end{center}

Word: $254314002544$,
\begin{center}
\begin{tikzcd}[column sep=2cm, row sep=0.1cm]
2^1 \arrow[r,  "f_{54314002544}"] & 2^{30}.
\end{tikzcd}
\end{center}

Word: $553032005532$,
\begin{center}
\begin{tikzcd}[column sep=2cm, row sep=0.1cm]
5 \arrow[r,  "f_{53032005532}"] & 2^{31}.
\end{tikzcd}
\end{center}

Word: $1550104015504$,
\begin{center}
\begin{tikzcd}[column sep=2cm, row sep=0.1cm]
1 \arrow[r,  "f_{550104015504}"]  & 2^{32}.
\end{tikzcd}
\end{center}

Word: $3540212035412$,
\begin{center}
\begin{tikzcd}[column sep=2cm, row sep=0.1cm]
3 \arrow[r,  "f_{540212035412}"] & 2^{33}.
\end{tikzcd}
\end{center}

Word: $11520424115224$,
\begin{center}
\begin{tikzcd}[column sep=3cm, row sep=0.1cm]
1\arrow[r, "f_{1520424115224}"] &2^{34}.
\end{tikzcd}
\end{center}

Word: $23441252234452$,
\begin{center}
\begin{tikzcd}[column sep=3cm, row sep=0.1cm]
2^1\arrow[r, "f_{3441252234452}"] & 2^{35}.
\end{tikzcd}
\end{center}

Word: $51322544513344$,
\begin{center}
\begin{tikzcd}[column sep=3cm, row sep=0.1cm]
5\arrow[r, "f_{1322544513344}"] & 2^{36}.
\end{tikzcd}
\end{center}

Word: $143045533431132$,
\begin{center}
\begin{tikzcd}[column sep=3cm, row sep=0.1cm]
1\arrow[r, "f_{43045533431132}"] & 2^{37}.
\end{tikzcd}
\end{center}

Word: $330135511302304$,
\begin{center}
\begin{tikzcd}[column sep=3cm, row sep=0.1cm]
3\arrow[r, "f_{30135511302304}"] & 2^{38}.
\end{tikzcd}
\end{center}

Word: $1100315423005012$,
\begin{center}
\begin{tikzcd}[column sep=3cm, row sep=0.1cm]
1\arrow[r, "f_{100315423005012}"] & 2^{39}.
\end{tikzcd}
\end{center}

Word: $2201035250014024$,
\begin{center}
\begin{tikzcd}[column sep=3cm, row sep=0.1cm]
2^1\arrow[r, "f_{201035250014024}"] & 2^{40}.
\end{tikzcd}
\end{center}

Word: $4402114540032052$,
\begin{center}
\begin{tikzcd}[column sep=3cm, row sep=0.1cm]
2^2\arrow[r, "f_{402114540032052}"] & 2^{41}.
\end{tikzcd}
\end{center}

Word: $13204233520104144$,
\begin{center}
\begin{tikzcd}[column sep=3cm, row sep=0.1cm]
1\arrow[r, "f_{3204233520104144}"] & 2^{42}.
\end{tikzcd}
\end{center}

Word: $30412511440212332$,
\begin{center}
\begin{tikzcd}[column sep=3cm, row sep=0.1cm]
3\arrow[r, "f_{0412511440212332}"] & 2^{43}.
\end{tikzcd}
\end{center}

Word: $101225423320425104$,
\begin{center}
\begin{tikzcd}[column sep=3cm, row sep=0.1cm]
1\arrow[r, "f_{01225423320425104}"] & 2^{44}.
\end{tikzcd}
\end{center}

Word: $202455251041254212$,
\begin{center}
\begin{tikzcd}[column sep=3cm, row sep=0.1cm]
2^1\arrow[r, "f_{02455251041254212}"] & 2^{45}.
\end{tikzcd}
\end{center}

Word: $405354542122552424$,
\begin{center}
\begin{tikzcd}[column sep=3cm, row sep=0.1cm]
2^2\arrow[r, "f_{05354542122552424}"] & 2^{46}.
\end{tikzcd}
\end{center}

Word: $1215153524245545252$,
\begin{center}
\begin{tikzcd}[column sep=3cm, row sep=0.1cm]
1\arrow[r, "f_2"]&2^3\arrow[r, "f_{15153524245545252}"] & 2^{47}.
\end{tikzcd}
\end{center}

Word: $2434351452535534544$,
\begin{center}
\begin{tikzcd}[column sep=3cm, row sep=0.1cm]
2^1\arrow[r, "f_4"]&2^4\arrow[r, "f_{34351452535534544}"] & 2^{48}.
\end{tikzcd}
\end{center}

Word: $5313143345515513532$,
\begin{center}
\begin{tikzcd}[column sep=3cm, row sep=0.1cm]
5\arrow[r, "f_{313143345515513532}"]&2^{49}.
\end{tikzcd}
\end{center}

Word: $15030331135435431504$,
\begin{center}
\begin{tikzcd}[column sep=3cm, row sep=0.1cm]
1\arrow[r, "f_{5030331135435431504}"]&2^{50}.
\end{tikzcd}
\end{center}

Word: $34101102315315303412$,
\begin{center}
\begin{tikzcd}[column sep=3cm, row sep=0.1cm]
3\arrow[r, "f_{4101102315315303412}"]&2^{51}.
\end{tikzcd}
\end{center}

Word: $112202205035035011224$,
\begin{center}
\begin{tikzcd}[column sep=3cm, row sep=0.1cm]
1\arrow[r, "f_{12202205035035011224}"]&2^{52}.
\end{tikzcd}
\end{center}

Word: $224404414114114022452$,
\begin{center}
\begin{tikzcd}[column sep=3cm, row sep=0.1cm]
2^1\arrow[r, "f_{24404414114114022452}"]&2^{53}.
\end{tikzcd}
\end{center}

Word: $453213232232232045344$,
\begin{center}
\begin{tikzcd}[column sep=3cm, row sep=0.1cm]
2^2\arrow[r, "f_{53213232232232045344}"]&2^{54}.
\end{tikzcd}
\end{center}

Word: $324324424525353241314531412$,
\begin{center}
\begin{tikzcd}[column sep=4cm, row sep=0.1cm]
3\arrow[r, "f_{24324424525353241314531412}"]&2^{69}.
\end{tikzcd}
\end{center}

Word: $4340341543153123340223221344$,
\begin{center}
\begin{tikzcd}[column sep=4cm, row sep=0.1cm]
2^2\arrow[r, "f_{340341543153123340223221344}"]&2^{72}.
\end{tikzcd}
\end{center}

Word: $543221154140454523304213143113125012$,
\begin{center}
\begin{tikzcd}[column sep=5cm, row sep=0.1cm]
5\arrow[r, "f_{43221154140454523304213143113125012}"]&2^{93}.
\end{tikzcd}
\end{center}

Word: $50420350225002300544352013210411505104$,
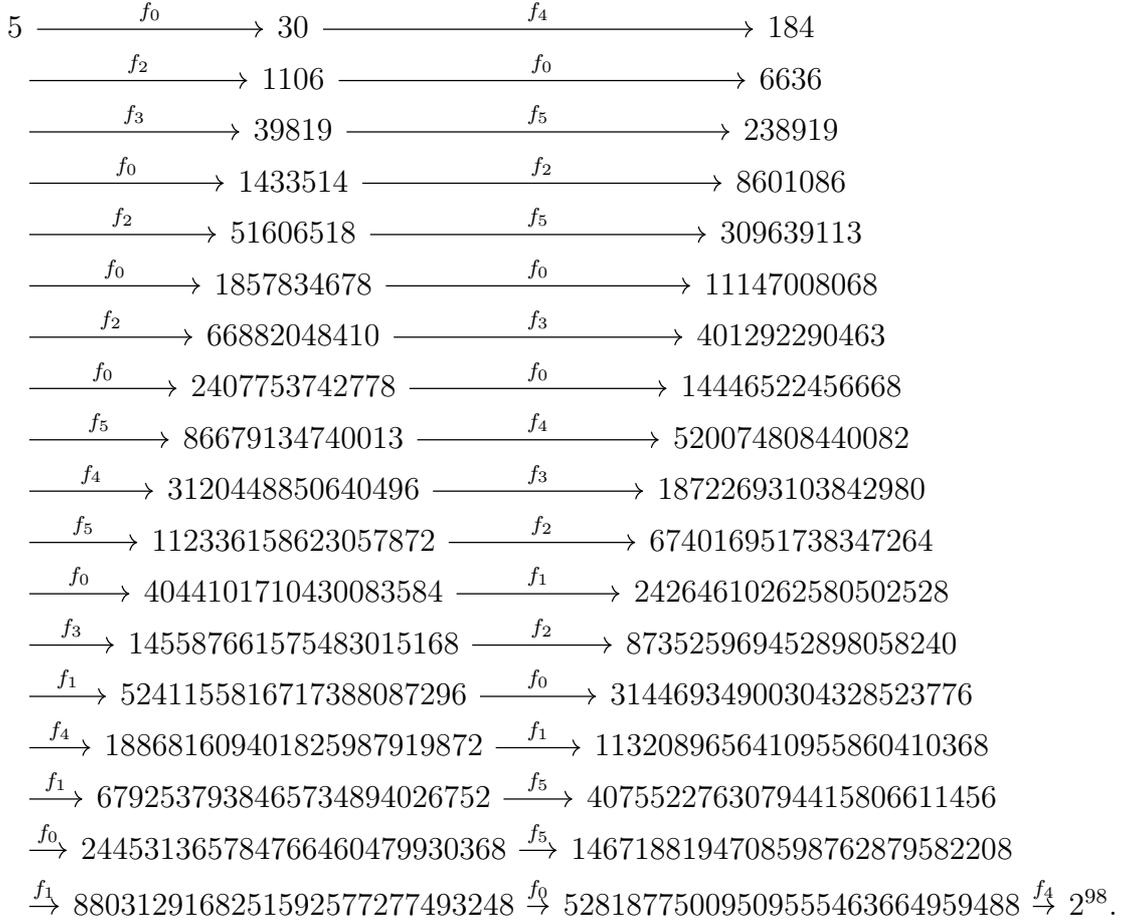
\begin{center}
\begin{tikzcd}[column sep=0.3cm, row sep=0.1cm]
5\arrow[r, "f_0"]&30\arrow[r, "f_4"]&184\\
\arrow[r, "f_2"]&1106\arrow[r, "f_0"]&6636\\
\arrow[r, "f_3"]&39819\arrow[r, "f_5"]&238919\\
\arrow[r, "f_0"]&1433514\arrow[r, "f_2"]&8601086\\
\arrow[r, "f_2"]&51606518\arrow[r, "f_5"]&309639113\\
\arrow[r, "f_0"]&1857834678\arrow[r, "f_0"]&11147008068\\
\arrow[r, "f_2"]&66882048410\arrow[r, "f_3"]&401292290463\\
\arrow[r, "f_0"]&2407753742778\arrow[r, "f_0"]&14446522456668\\
\arrow[r, "f_5"]&86679134740013\arrow[r, "f_4"]&520074808440082\\
\arrow[r, "f_4"]&3120448850640496\arrow[r, "f_3"]&18722693103842980\\
\arrow[r, "f_5"]&112336158623057872\arrow[r, "f_2"]&674016951738347264\\
\arrow[r, "f_0"]&4044101710430083584\arrow[r, "f_1"]&24264610262580502528\\
\arrow[r, "f_3"]&145587661575483015168\arrow[r, "f_2"]&873525969452898058240\\
\arrow[r, "f_1"]&5241155816717388087296\arrow[r, "f_0"]&31446934900304328523776\\
\arrow[r, "f_4"]&188681609401825987919872\arrow[r, "f_1"]&1132089656410955860410368\\
\arrow[r, "f_1"]&6792537938465734894026752\arrow[r, "f_5"]&40755227630794415806611456\\
\arrow[r, "f_0"]&244531365784766460479930368\arrow[r, "f_5"]&1467188194708598762879582208\\
\arrow[r, "f_1"]&8803129168251592577277493248\arrow[r, "f_0"]&52818775009509555463664959488\arrow[r, "f_4"]&2^{98}.
\end{tikzcd}
\end{center}
\end{exa}

One can see from Example \ref{exa2} the following conclusions on $\sim$ and $\stackrel{par}{\sim}$. 
\begin{itemize}
\item $2^3, 2^6, 2^8, 2^{11}, 2^{13}, 2^{16}, 2^{19}, 2^{21}, 2^{24},2^{26}, 2^{29}, 2^{32}, 2^{34},2^{37}, 2^{39}, 2^{42}, 2^{44}, 2^{47}, 2^{50}, 2^{52}$ are all descendants of $2^0$. Among others  $2^3, 2^6, 2^8, 2^{11}, 2^{13}, 2^{19}, 2^{21}, 2^{24},2^{26}, 2^{32}, 2^{34},2^{37}, 2^{39},\\ 2^{42}, 2^{44}, 2^{50}, 2^{52}$ are all children, while $2^{16}, 2^{29}, 2^{47}$ are only descendants instead of children.

\item $2^4, 2^9, 2^{14}, 2^{17}, 2^{22}, 2^{27}, 2^{30}, 2^{35}, 2^{40}, 2^{45}, 2^{48}, 2^{53}$ are all descendants of $2^1$. Among others  $2^4, 2^9, 2^{14}, 2^{22}, 2^{27}, 2^{30}, 2^{35}, 2^{40}, 2^{45}, 2^{53}$ are all children, while $2^{17}, 2^{48}$ are only descendants instead of children.

\item $2^{10}, 2^{15}, 2^{23}, 2^{28}, 2^{41}, 2^{46}, 2^{54}, 2^{72}$ are all descendants of $2^2$. All are in fact also children of it. 

\item $2^{16}, 2^{29}, 2^{47}$ are all descendants of $2^3$. All are in fact also children of it. 

\item $2^{17}, 2^{48}$ are both descendants of $2^4$. Both are in fact also children of it.
\end{itemize}

Counting chains of maximal length, there are 5 chains 
$$2^0\stackrel{par}{\sim}2^3\stackrel{par}{\sim}2^{16},$$
$$2^0\stackrel{par}{\sim}2^3\stackrel{par}{\sim}2^{29},$$
$$2^0\stackrel{par}{\sim}2^3\stackrel{par}{\sim}2^{47},$$
$$2^1\stackrel{par}{\sim}2^4\stackrel{par}{\sim}2^{17},$$
$$2^1\stackrel{par}{\sim}2^4\stackrel{par}{\sim}2^{48},$$ 
of length $3$, 21 chains 
$$2^1\stackrel{par}{\sim}2^9, \cdots, 2^2\stackrel{par}{\sim}2^{10}, \cdots, 2^3\stackrel{par}{\sim}2^{47}$$
of length $2$ and 32 chains 
$$2^5, \cdots, 2^{98}$$
of length $1$ exhibited in Example \ref{exa2}.

One can see from Example \ref{exa2} that every length-2 word is the prefix of at least one word corresponding to some (finite) chain in Example \ref{exa2}. Note that length-2 word beginning with $0$ can be fabricated from the exhibited words in the example. For example, if we add $0$ ahead of the word $52$ (corresponding to the chain $2^5$ of length $1$), we get $052$ (corresponding to the same chain whose induced numbers reside in the depth-2 fundamental interval $I_{05}$).   

The next example should be understood along with  Proposition \ref{thm2}.

\begin{exa}\label{exa3}
Let $S=\{2,5,7,11\}$ and $b=3$. We consider the finite factorial number $\xi_{S,b,\bold{a}_0}^{FF,-}$ with respect to the word $\bold{a}_0=1$ here. We provide the numerical values of the nodes $\{m_i\}_{i=1}^5$, their exponents $\{\bold{r}_i\}_{i=1}^5$ together with the successive words $\{\bold{a}_i\}_{i=1}^5$ as in the proof of Proposition \ref{thm2} with respect to the prescribed $S, b, \bold{a}_0$. We accompany  $m_1$ with the second least child $m_{1,<}$ (which is the least one larger than $m_1$) of the word $\bold{a}_0$, its exponent $\bold{r}_{1,<}$ and the corresponding word $\bold{a}_{1,<}$. For every $1\leq i\leq 4$, we accompany the above terms with the second least child $m_{i+1,<}$ of $m_i$, its exponent $\bold{r}_{i+1}$ and the corresponding word $\bold{a}_{i+1}$.     

Index: i=1,
\begin{center}
$m_1=4$, $\bold{r}_1=(2,0,0,0)$, $\bold{a}_1=1$.

$m_{1,<}=5$, $\bold{r}_{1,<}=(0,1,0,0)$, $\bold{a}_{1,<}=2$.
\end{center}

Index: i=2,
\begin{center}
$m_2=14$, $\bold{r}_2=(1,0,1,0)$, $\bold{a}_2=2$.

$m_{2,<}=40$, $\bold{r}_{2,<}=(3,1,0,0)$, $\bold{a}_{2,<}=2$.
\end{center} 

Index: i=3,
\begin{center}
$m_3=44$, $\bold{r}_3=(2,0,0,1)$, $\bold{a}_3=2$.

$m_{3,<}=128$, $\bold{r}_{3,<}=(7,0,0,0)$, $\bold{a}_{3,<}=02$.
\end{center}

Index: i=4,
\begin{center}
$m_4=400$, $\bold{r}_4=(4,2,0,0)$, $\bold{a}_4=11$.

$m_{4,<}=1210$, $\bold{r}_{4,<}=(1,1,0,2)$, $\bold{a}_{4,<}=211$.
\end{center}

Index: i=5,
\begin{center}
$m_5=875000$, $\bold{r}_5=(3,6,1,0)$, $\bold{a}_5=0021102$.

$m_{5,<}=7884800$, $\bold{r}_{5,<}=(12,2,1,1)$, $\bold{a}_{5,<}=120220122$.
\end{center}
\end{exa}  

One can see that 
$$\xi_{S,b,\bold{a}_0}^{FF,-}=(.1122110021102\cdots)_b$$
for $S=\{2,5,7,11\}, b=3$ and $\bold{a}_0=1$ from Example \ref{exa3}.

\section*{\raggedright{Appendix}}

There are various generalizations of the Schmidt Subspace Theorem as well as their applications in different contexts, while the following version suffices for our aim in this work (refer to Schlickewei \cite{Schl}).
\begin{Schmidt Subspace Theorem}
Let $S\subset P$ be a finite set of primes and $d\in \mathbb{N}$ be some integer at least $2$. Let $\{\mathcal{L}_{A,i}(\bold{x})\}_{1\leq i\leq d}$ be a sequence of independent linear forms with algebraic coefficients and $\{\mathcal{L}_{R,i}(\bold{x})\}_{1\leq i\leq d}$ be a sequence of independent linear forms with rational coefficients for $\bold{x}=(x_1,x_2,\cdots,x_d)$. The for any $\epsilon_{1.2}>0$, the integral solutions $\bold{x}\in\mathbb{Z}^d$ of the inequality
$$\prod_{1\leq i\leq d}|\mathcal{L}_{A,i}(\bold{x})|\prod_{1\leq i\leq d, s\in S}|\mathcal{L}_{R,i}(\bold{x})|_s\leq |\bold{x}|^{-\epsilon_{1.2}}$$
lie in finitely many proper subspaces of $\mathbb{Q}^d$.
\end{Schmidt Subspace Theorem}

To state the Evertse's Theorem (refer to \cite{Eve}) on the bounds of $S$-unit equations on an algebraic number field $K$, we introduce some notions first.  An equivalent class of non-trivial valuations on $K$ is termed a \emph{place} (which is finite if there are only non-archimedean valuations and infinite otherwise). Let $S$ be a finite set of places containing all the infinite ones (there are only finitely many infinite places on $K$). The $S$-\emph{units} are collections of numbers such that if some valuation takes $1$ at the number then it does not belong to $S$.  
\begin{Evertse's Theorem}
Let $K$ be an algebraic number field of degree $d$ and $\{\alpha,\beta,\gamma\}\subset K\setminus\{0\}$, let $S$ be a finite set of places. Then the equation 
$$\alpha x+\beta y=\gamma$$
with the unknows being $S$-units admits at most $7^{d+2\#S}$ solutions.
\end{Evertse's Theorem}

Although the bound seems to be a big one in case $d$ or $\#S$ being large, for most equations it can be decreased to $2$ in fact, refer to Evertse-Gy\"ory-Stewart-Tijdeman \cite{EGST}. For the upper bounds on the nondegenerate solutions of multidiemntional $S$-unit equations one can refer to Evertse-Schlickewei-Schmidt \cite{ESS}.

\begin{Bird's Theorem}
Let $N=\{n_i\}_{i=1}^\infty\subset \mathbb{N}_+$ and $b\geq 2$ be some base. Assume $\lim_{i\rightarrow\infty} \cfrac{n_{i+1}}{n_i}=\alpha$ exists (could be $\infty$). Then $N$ is extendable in base $b$ if either of the following two conditions is satisfied,
\begin{enumerate}[(1).]
\item $\alpha=1$ or $b\notin \{\alpha^{j/k}\}_{j,k\in\mathbb{Z}\setminus\{0\}}$,

\item $\alpha=\infty$ and $\lim_{i\rightarrow\infty} \cfrac{n_in_{i+2}}{n_{i+1}^2}=\alpha$.
\end{enumerate}
\end{Bird's Theorem}

Both conditions imply the set $\Big\{(\cfrac{\log n_i}{\log b})\Big\}_{i=1}^\infty$ is dense in the unit interval $(0,1)$, through an argument of Cassels (\cite{Cas}).


\begin{thebibliography}{88}

\bibitem[AB1]{AB1} B. Adamczewski and Y. Bugeaud, On the complexity of algebraic numbers. I. Expansions in integer bases, Ann. Math.,  165 (2),  547-565, 2007.

\bibitem[AB2]{AB2} B. Adamczewski and Y. Bugeaud, Palindromic continued fractions, Ann. Inst. Fourier (Grenoble), 57,  1557-1574, 2007.

\bibitem[AB3]{AB3} B. Adamczewski and Y. Bugeaud, Real and $p$-adic expansions involving symmetric patterns, Inter. Math. Res. Not., 2006, 1-17, 2006.

\bibitem[ABL]{ABL} B. Adamczewski, Y. Bugeaud and F. Luca, Sur la complexit\'e des nombres alg\'ebriques, C. R. Acad. Sci. Paris, 339, 11-14, 2004.

\bibitem[ADM]{ADM} B. Adamczewski, M. Drmota and C. M\"ullner, (Logatihmic) densities for automatic sequences along prime and squares, Trans. Amer. Math. Soc., 375 (1), 455-499, 2022.

\bibitem[AF]{AF} B. Adamczewski and C. Faverjon,	Mahler's method in several variables and finite automata, arXiv:2012.08283 [math.NT], 2020.

\bibitem[AS1]{AS1} J.-P. Allouche and J. Shallit, Automatic Sequences: Theory, Applications, Generalizations, Cambridge University Press, 2003.

\bibitem[AS2]{AS2} J.-P. Allouche and J. Shallit, Complexit\'e des suites de Rudin-Shapiro g\'en\'eralis\'ees, J. Th\'eor. Nombres Bordeaux, 5 (2),  283-302, 1993.

\bibitem[AZ]{AZ} J.-P. Allouche and L. Zamboni, Algebraic irrational binary numbers cannot be fixed points of non-trivial constant length or primitive morphisms, J. Number Theory, 69 (1), 119-124, 1998.

\bibitem[BC]{BC} K. Boklan and J. Conway, Expect at Most One Billionth of a New Fermat Prime!,
The Mathematical Intelligencer, 39, 3-5, 2017.

\bibitem[Bel1]{Bel1} J. Bell, Logarithmic frequency in morphic sequences, J. Th\'eor. Nombres Bordeaux, 20 (2), 227-241, 2008. 

\bibitem[Bel2]{Bel2} J. Bell, The upper density of an automatic set is rational, J. Th\'eor. Nombres Bordeaux, 32 (2), 585-604, 2020.

\bibitem[BHZ]{BHZ} V. Berth\'e, C. Holton, and L. Zamboni, Initial powers of Sturmian sequences, Acta Arith., 122,  315-347, 2006.

\bibitem[Bir]{Bir} R. Bird, Integers with given initial digits, Amer. Math. Monthly, 79, 367-370, 1972.

\bibitem[Bor]{Bor} \`E. Borel, Sur les chiffres d\'ecimaux de $\sqrt{2}$ et divers probl\`emes de probabilit\'es en cha\^ine, C. R. Acad. Sci. Paris, 230,  591-593, 1950.

\bibitem[Cas]{Cas} J. Cassels, An Introduction to Diophantine Approximation, Cambridge Tracts in Math., No. 45, C. U. P., England, 1965.

\bibitem[Cob]{Cob} A. Cobham, Uniform tag sequences, Math. Syst. Theory, 6, 164-192, 1972.

\bibitem[DMR]{DMR} M. Drmota, C. Mauduit and J. Rivat, The sum-of-digits function of polynomial sequences, J. Lond. Math. Soc.,  84 (2), 81-102, 2011.

\bibitem[DS]{DS} H. Davenport and W. Schmidt, Approximation to real numbers by algebraic integers, Acta Arith., 15, 393-416, 1969.

\bibitem[Dys]{Dys} F. Dyson, The approximation to algebraic numbers by rationals, Acta Math., 79, 225-240, 1947.

\bibitem[EGST]{EGST} J.-H. Evertse, K. Gy\"ory, C. Stewart and R. Tijdeman, On $S$-unit equations in two unknowns, Invent. Math., 92 (3), 461-478, 1988.

\bibitem[ESS]{ESS} J.-H. Evertse, H. Schlickewei and W. Schmidt, Linear equations in variables which lie in a multiplicative group, Ann. Math., 155 (3),  807-836, 2002.

\bibitem[ET]{ET} P. Erd\"os and P. Tur\'an, On a problem in the elementary theory of numbers, Amer. Math. Monthly, 41, 608-611, 1934.

\bibitem[Eve]{Eve} J.-H. Evertse, On equations in S-units and the Thue-Mahler equation, Invent. Math., 75, 561-584, 1984.

\bibitem[Fis1]{Fis1} S. Fischler, Palindromic prefixes and episturmian words, J. Combin. Theory. Ser. A, 113 (7),  1281-1304, 2006.

\bibitem[Fis2]{Fis2} S. Fischler, Palindromic Prefixes and Diophantine Approximation, Monatsh. Math., 151, 11-37, 2007.

\bibitem[FM]{FM} S. Ferenczi and C. Mauduit, Transcendence of numbers with a low complexity expansion, J. Number Theory, 67 (2),  146-161, 1997.

\bibitem[Gel]{Gel} A. Gelfond, Sur les nombres qui ont des propri\'et\'es additives et multiplicatives donn\'ees, Acta Arith., 13, 259-265, 1967.

\bibitem[Goo]{Goo}  I. Good, The fractional dimensional theory of continued fractions, Math. Proc. Camb. Phil. Soc., 37, 199-228, 1941.

\bibitem[GST]{GST} K. Gy\"ory, C. Stewart and R. Tijdeman, On prime factors of sums of integers I, Compositio Math., 59 (1),  81-88, 1986.

\bibitem[HW]{HW} G. Hardy and E. Wright, An Introduction to the Theory of Numbers, Oxford University Press, 5th ed., 1985.

\bibitem[Iwa]{Iwa} H. Iwaniec, Almost-primes represented by quadratic polynomials, Invent. Math., 47 (2), 178-188, 1978.

\bibitem[KN]{KN} L. Kuipers and H. Niederreiter, Uniform distribution of sequences, Wiley-Interscience, John
Wiley \& Sons, New York, 1974. 

\bibitem[Lio]{Lio} J. Liouville,  Sur des classes tr\`es \'etendues de quantit\'es dont la valeur n’est ni alg\'ebrique, ni m\^eme reductible \`a des irrationelles alg\'ebriques, J. Math. Pur. App., 1 (16), 133-142, 1851.

\bibitem[LM]{LM} D. Lewis and K. Mahler, Representation of integers by binary forms, Acta Arith.,  6, 333-363, 1961.

\bibitem[LV1]{LV1} J. Loxton and A. van der Poorten, Arithmetic properties of the solutions of a class of functional equations, J. Reine Angew. Math., 330, 159-172, 1982.

\bibitem[LV2]{LV2} J. Loxton and A. van der Poorten, Arithmetic properties of automata: regular sequences, J. Reine Angew. Math., 392, 57-69, 1988.

\bibitem[Mah]{Mah} K. Mahler,  Arithmetische Eigenschaften der L\"osungen einer Klasse von Funktionalgleichungen, Math. Annalen, 101, 342-366, 1929, Corrigendum, 103, 532, 1930.

\bibitem[MR1]{MR1} C. Mauduit and J. Rivat, La somme des chiffres des carr\'es, Acta Math., 203 (1), 107-148, 2009.

\bibitem[MR2]{MR2} C. Mauduit and J. Rivat, Sur un probl\`eme de Gelfond: la somme des chiffres des nombres premiers, Ann. Math., 171 (3), 1591-1646, 2010.


\bibitem[Schl]{Schl} H. Schlickewei, The p-adic Thue-Siegel-Roth-Schmidt theorem, Arch. Math., XXIX, 267-270, 1977.

\bibitem[Schm]{Schm} W. Schmidt, Norm form equations, Ann. Math.,  96 (3), 526-551, 1972.

\bibitem[Rot]{Rot} K. Roth, Rational approximations to algebraic numbers, Mathematika, 2 (1), 1-20, 1955, Corrigendum, 168, 1955.

\bibitem[RW1]{RW1} D. Roy and M. Waldschmidt, Approximation diophantienne et ind\'ependance alg\'ebrique de logarithmes, Ann. Sci. \'Ecole Norm. Sup., 30 (6), 753-796, 1997.

\bibitem[RW2]{RW2} D. Roy and M. Waldschmidt, Diophantine approximation by conjugate algebraic integers, Compositio Math., 140, 593-612, 2004.

\bibitem[Sie]{Sie} C. Siegel, Approximation algebraischer Zahlen, Math. Zei., 10, 173-213, 1921.

\bibitem[Thu]{Thu} A. Thue, \"Uber Annäherungswerte algebraischer Zahlen, J. Reine Angew. Math., 135, 284-305, 1909. 

\bibitem[Wey]{Wey} H. Weyl, \"Uber die Gleichverteilung von Zahlen mod. Eins, Math. Annalen, 77, 13-352, 1916.

\end{thebibliography}
\end{document}